\documentclass[11pt,leqno]{article}

\usepackage{amsmath,amsfonts,amscd,amssymb,theorem}

\long\def\comment#1\endcomment{}

\makeatletter
\begingroup
\gdef\th@dotted{\normalfont\itshape
  \def\@begintheorem##1##2{%
        \item[\hskip\labelsep \theorem@headerfont ##1\ ##2.]}%
\def\@opargbegintheorem##1##2##3{%
   \item[\hskip\labelsep \theorem@headerfont ##1\ ##2\ (##3).]}}
\endgroup
\makeatother

\theoremstyle{dotted}

\newtheorem{theorem}{Theorem}[section]
\newtheorem{lemma}[theorem]{Lemma}

\newtheorem{prop}[theorem]{Proposition}

\makeatletter
\begingroup
\gdef\th@upshape{\normalfont
  \def\@begintheorem##1##2{%
        \item[\hskip\labelsep \theorem@headerfont ##1\ ##2.]}%
\def\@opargbegintheorem##1##2##3{%
   \item[\hskip\labelsep \theorem@headerfont ##1\ ##2\ (##3).]}}
\endgroup
\makeatother

\theoremstyle{upshape}

\newtheorem{defn}[theorem]{Definition}
\newtheorem{remark}[theorem]{Remark}

\makeatletter
\renewcommand{\subsection}{\@startsection{subsection}{2}{0pt}{-3ex
plus -1ex minus -0.2ex}{-2mm plus -0pt minus
-2pt}{\normalfont\bfseries}} \makeatother

\makeatletter
\@addtoreset{equation}{section}
\makeatother

\newcommand{\cntrct}                
{\hspace{2pt}\raisebox{1pt}{\text{$\lrcorner$}}\hspace{2pt}}

\newcommand{\proof}[1][Proof.]{\smallskip\noindent{\em #1}}
\def\endproof{\hfill\ensuremath{\square}\par\medskip}

\renewcommand{\labelenumi}{{\normalfont(\roman{enumi})}}

\def\eqref#1{\thetag{\ref{#1}}}

\let\latexref=\ref
\def\ref#1{{\normalfont{\latexref{#1}}}}

\newcommand{\wt}{\widetilde}
\newcommand{\wh}{\widehat}

\newcommand{\onto}{\twoheadrightarrow} 

\newcommand{\hdot}{{\:\raisebox{3pt}{\text{\circle*{1.5}}}}}
\newcommand{\C}{{\mathbb C}}

\newcommand{\Z}{{\mathbb Z}}

\newcommand{\A}{{\cal A}}
\newcommand{\CC}{{\cal C}}
\newcommand{\D}{{\cal D}}
\newcommand{\J}{{\cal J}}

\newcommand{\M}{{\cal M}}
\newcommand{\Mb}{\overline{\M}}

\newcommand{\T}{{\cal T}}

\newcommand{\U}{{\cal U}}
\newcommand{\V}{{\cal V}}
\newcommand{\W}{{\cal W}}

\newcommand{\X}{{\mathfrak X}}

\newcommand{\E}{{\cal E}}

\newcommand{\G}{{\sf G}}
\newcommand{\Gb}{\overline{\G}}
\newcommand{\Gt}{\widetilde{\G}}
\newcommand{\HH}{{\sf H}}
\newcommand{\WW}{{\sf W}}

\newcommand{\calo}{{\cal O}}

\newcommand{\m}{{\mathfrak m}}
\newcommand{\g}{{\mathfrak g}}

\newcommand{\h}{{\mathfrak h}}

\newcommand{\hh}{{\cal H}}
\newcommand{\gl}{{\mathfrak{g}\mathfrak{l}}}
\newcommand{\gm}{{\mathbb{G}_m}}

\renewcommand{\phi}{\varphi}

\renewcommand{\dim}{\operatorname{\sf dim}}

\newcommand{\id}{\operatorname{\sf id}}

\newcommand{\Hom}{\operatorname{Hom}}
\newcommand{\Ext}{\operatorname{Ext}}
\newcommand{\ext}{\operatorname{{\cal E}{\it xt}}}
\newcommand{\spl}{\operatorname{{\cal S}{\it pl}}}

\newcommand{\Aut}{\operatorname{{\sf Aut}}}
\newcommand{\Sympl}{\operatorname{{\sf Symp}}}

\newcommand{\Der}{\operatorname{Der}}

\newcommand{\Ker}{\operatorname{{\sf Ker}}}
\newcommand{\Coker}{\operatorname{{\sf Coker}}}

\newcommand{\Spec}{\operatorname{Spec}}
\newcommand{\Spf}{\operatorname{Spf}}

\newcommand{\Per}{\operatorname{{\sf Per}}}

\newcommand{\loc}{\operatorname{{\sf Loc}}}
\newcommand{\Lin}{\operatorname{{\sf Lin}}}

\newcommand{\whotimes}{\operatorname{\wh{\otimes}}}

\newcommand{\GH}{{\langle G,\h \rangle}}
\newcommand{\GHh}{{\langle G_1,\h_1 \rangle}}
\newcommand{\GHhh}{{\langle G'_1,\h_1 \rangle}}
\newcommand{\GHt}{{\langle G_0,\h_0 \rangle}}
\newcommand{\vv}{{\langle V,V \rangle}}
\newcommand{\WA}{{\langle \Aut \A,\WW \rangle}}
\newcommand{\WSA}{{\langle \Sympl \A,\HH \rangle}}
\newcommand{\WDA}{{\langle \Aut D, \Der D \rangle}}
\newcommand{\WDAp}{{\langle (\Aut D)_p, (\Der D)_p \rangle}}
\newcommand{\WDApt}{{\langle (\Aut D)_{p+1}, (\Der D)_{p+1} \rangle}}

\newcommand{\dg}{\dagger}

\title{Fedosov quantization in algebraic context}

\author{R. Bezrukavnikov\thanks{Partially supported by NSF grant
DMS-0071967} and D. Kaledin\thanks{Partially supported by CRDF grant
RM1-2354-MO02.}}

\date{{\em To B.L. Feigin on his 50th anniversary}}

\begin{document}

\maketitle

\begin{abstract} 
We consider the problem of quantization of smooth symplectic
varieties in the algebro-geometric setting. We show that, under
appropriate cohomological assumptions, the Fedosov quantization
procedure goes through with minimal changes. The assumptions are
satisfied, for example, for affine and for projective varieties. We
also give a classification of all possible quantizations.
\end{abstract}

\tableofcontents

\section*{Introduction}

Let $V$ be a finite-dimensional vector space equiped with a
non-degenerate $2$-form. The algebra $S(V)$ of polynomial functions
on $V$ admits a well-known non-commutative one-parameter deformation
$S(V)[[h]]$ called the Weyl algebra. The problem of deformation
quantization consist in generalizing this construction to a
deformation of the sheaf $\A(M)$ of functions on an arbitrary smooth
symplectic manifold $M$. More precisely, one wants to know whether
there exists a deformation with prescribed properties, and how many
such deformations there are.

The reader will immediately notice that our language is ambiguous:
``smooth manifold'' can mean either a $C^\infty$-manifold, or a
holomorphic manifold, or a smooth algebraic variety -- over $\C$ or
over some other field, possibly of positive characteristic. This is
intentional: the problem of deformation quantization makes perfect
sense in all these situations.

When the problem was posed several decades ago, it soon became clear
that the standard deformation theory methods take one only so
far. General nonsense gives a series of obstruction classes lying in
a certain group. However, this group is usually non-trivial. Thus to
quantize a manifold, it is necessary to actually prove that the
obstruction classes themselves vanish.

After a hiatus of several years, the problem was finally completely
solved in the early 1980-ies independently by M. De Wilde-P. Lecomte
and by B. Fedosov (see \cite{DWL}, \cite{fedo} and a classic
exposition of these results by P. Deligne in \cite{D}). The answer
is that a quantization always exists, and that the space of all
quantizations admits a simple description.

Both De Wilde-Lecomte and Fedosov worked with $C^\infty$-manifolds,
by $C^\infty$ methods. So did Deligne. When one looks at the proofs,
though, one is tempted to think that the $C^\infty$ context is not
really essential -- one only needs the vanishing of certain
cohomology groups. This is implicit in \cite{DWL} and \cite{fedo},
and less implicit in the gerb-theoretic version of the proof given
in \cite{D}. However, Deligne does not state the necessary
cohomology vanishing conditions either. Instead, he uses the
softness of certain non-abelian group sheaves. Thus one cannot
directly generalize either of the existing proofs to the holomorphic
or algebraic setting -- while there is a strong feeling that the
results themselves should hold.

In the eight years which passed since the publication of \cite{D},
the deformation quantization has been much better understood, and
now there seems to be no doubt among experts as to what happens in
the holomorphic and in the algebraic setting (at least in
characteristic $0$). Some proofs are actually published. In
particular, R. Nest and B. Tsygan have given in \cite{NT} a complete
proof in the holomorphic case. They have also specified the
cohomology vanishing condition which one needs to impose on the
manifold in order for the argument to work.

\bigskip

However, it seems that the algebraic case still remains a folk
knowledge, with no references in the literature. Thus a write-up of
a purely algebraic proof would be useful. This is what the present
paper is intended to be.

\bigskip

The results in the paper were discovered while trying to apply
deformation quantization to a concrete algebro-geometric
problem. The authors are definitely not experts in the field, and we
lay no claim whatsoever to the novelty of our results. Moreover,
even our approach is essentially the same as Fedosov's, although
retold in a more algebraic language. The main technical tool is the
bundle of formal coordinate systems and the associated bundle of
jets (``formal geometry'' in the language of I.M. Gel'fand). We are
deeply gratefull to B. Feigin who suggested this approach to us and
more or less explained what to do.

We would like to mention also a recent paper \cite{KV}, by the
second author jointly with M. Verbitsky, which contains certain
results on purely commutative deformation of symplectic manifolds --
more or less, a generalization of the unobstructedness theorem of
F. Bogomolov \cite{bogo}. The methods used there are different and
somewhat simpler. However, the final result is completely parallel
to what one has for quantizations. In particular, the required
cohomology vanishing is precisely the same. In the latter part of
the paper, we explain this similarity and show how to join
quantizations and symplectic deformation into a single
partially-commutative deformation of the manifold in question.

Finally -- our proof only works in characteristic $0$. What happens
in characteristic $p > 0$? There are important reasons to study this
question, and we believe that it is possible to prove some sort of a
general statement. However, if one wants to apply our methods, one
has to modify them in quite an essential way. We plan to return to
this in future research.

\subsection*{Acknowledgements.} Aside from the very helpful, indeed
crucial contribution of B.L. Feigin which we have already mentioned,
we would like to thank A. Beilinson, F. Bogomolov, M. Finkelberg,
V. Ginzburg, I. Mirkovich, D. Tamarkin, B. Tsygan and M. Verbitsky
for stimulating discussions and suggestions. Part of the work was
done during the second author's visit to University of Chicago and
to Northwestern University. The hospitality of these institutions is
gratefully acknowledged. Another, and larger by far, part of the
work was done during the first author's sojourn in Moscow while
waiting for a US visa. We would like to extend our heartfelt thanks
to the State Department and other government agencies of the USA for
promoting research by obstructing travel.

\subsection*{Note added in proof.} A short time after the first
version of this paper was posted to the web, there appeared a
preprint \cite{Y} by A. Yekutieli devoted to a related
problem. A. Yekutieli also read our paper and indicated some
omissions which we now fix. Additonally, we would like to thank the
referee for the detailed report and important suggestions, and the
people at the Moscow Math Journal for their remarkable patience with
our extensive last-minute rewrites.

\section{Statements and definitions.}

\subsection{Notation.}\label{dfn}
Fix once and for all a base scheme $S$. Throughout the paper we will
assume that $S$ is a scheme of finite type over a fixed field $k$ of
charasteristic $0$. The most important case for us is $S = \Spec k$,
a point. However, all the proofs work for non-trivial schemes just
as well, and in some applications it is convenient to have the
results available in a more general setting.

By an {\em $S$-manifold} $X$ we will understand a scheme $X/S$ of
finite type and smooth over $S$ -- that is, we require that $X$ is
flat over $S$ and the relative cotangent sheaf $\Omega^1_{X/S}$ is a
locally free coherent sheaf. By the {\em dimension} of an
$S$-manifold we will understand the relative dimension $X/S$, which
coincides with the rank of the flat sheaf $\Omega^1_{X/S}$. For an
$S$-manifold $X$, one defines the (relative) de Rham complex
$\Omega^\hdot_{X/S}$ and its hypercohomology, known as the {\em de
Rham cohomology groups} $H^\hdot_{DR}(X/S)$. When $S = \Spec \C$ is
the complex point, the de Rham cohomology groups are known to
coincide with the topological cohomology groups $H^\hdot(X,\C)$.

By a ``vector bundle'' we will understand a ``locally free coherent
sheaf''.  For a vector bundle $\E$ on $X$, we define a (relative)
{\em flat connection} $\nabla$ on $\E$ as a differential operator
$\nabla:\E \to \E \otimes \Omega^1(X/S)$ satisfying the usual
compatibilities. Flat vector bundles on $X$ form a tensor abelian
category, with unit object $\calo_X$ (with the tautological
connection). For every flat vector bundle $\E$, one defines the
(relative) de Rham cohomology groups $H^\hdot_{DR}(X,\E)$. The de
Rham cohomology of the unit bundle $\calo_X$ coincide with the de
Rham cohomology groups $H^\hdot_{DR}(X)$.

The coherent cohomology groups $H^\hdot(X,\E)$ of a vector bundle
$\E$ on $X$ can be interpreted as relative de Rham cohomology by
means of relative jet bundles. To define those, let $\wh{\Delta}$ be
the completion of the fibered product $X \times_S X$ along the
diagonal $\Delta \subset X \times_S X$, and let
$\pi_1,\pi_2:\wh{\Delta} \to X$ be the projections onto the first
and the second factor. Then the jet bundle $J^\infty\E$ is given by
$$
J^\infty\E = \pi_{1*}\pi_2^*\E.
$$
The jet bundle carries a natural flat connection. The sheaf of its
flat sections coincide with the sheaf $\E$, and the de Rham
cohomology $H^\hdot_{DR}(X,J^\infty\E)$ is canonically isomorphic to
$H^\hdot(X,\E)$.

Note that a jet bundle $J^\infty\E$ is not finitely generated as a
sheaf of $\calo_X$-modules, thus not coherent. To be able to work
with jet bundles, we have to complete the category of coherent
sheaves on $X$ by adding countable projective limits. The resulting
category of pro-coherent sheaves is a tensor abelian category
(although it no longer has good duality properties). For the details
of the completion procedure, see \cite{D2}. As an additional bonus
for working with the completed category, we can interpret the de
Rham cohomology groups $H^\hdot_{DR}(X,\E)$ of a flat vector bundle
$\E$ as the $\Ext^\hdot$-groups from $\calo_X$ to $\E$ (in the usual
category, this is not necessarily true even for $\E \cong \calo_X$).
For the proof, it suffices to consider the de Rham type resolution
of $\E$ by jet bundles $J^\infty\Omega^\hdot_X \otimes \E$.

To simplify notation, we will often drop $S$ from the formulas and
omit the word ``relative'' in the statements. The reader should
always keep in mind that everything on $X$ is understood relatively
over $S$. Moreover, we will drop the prefix ``pro'' whenever there
is no danger of confusion.

\subsection{Assumptions.}
Let $X$ be an $S$-manifold. All our results will be valid under the
following assumption.

\begin{defn}\label{adm} 
The manifold $X$ is called {\em admissible} if the canonical map
$$
H_{DR}^i(X) \to H^i(X,\calo_X)
$$
from the de Rham cohomology $H^i_{DR}(X)$ to the cohomology
$H^i(X,\calo_X)$ of the structure sheaf $\calo_X$ is surjective for
$i=1,2$.
\end{defn}

In the case when $S = \Spec k$ is a point, examples of admissible
manifolds are:
\begin{enumerate}
\renewcommand{\labelenumi}{(\alph{enumi})}
\item A projective manifold $X$ -- admissibility follows from the
  Hodge theory.
\item A smooth projective resolution $X \to Y$ of a singular affine
  variety $Y$ such that the canonical bundle $K_X$ is trivial -- we
  have $H^i(X,\calo_X) = 0$ for all $i \geq 1$ by the
  Grauert-Riemenschneider Vanishing Theorem.
\end{enumerate}
For an $S$-manifold $X$, we will denote by $H_F^\hdot(X)$ the
hypercohomology of the first piece $F^1\Omega^\hdot_X$ of the de
Rham complex $\Omega_X^\hdot$ with respect to the {\em filtration
b\^{e}te} -- in other words, the third term in the natural
cohomology long exact sequence
$$
\begin{CD}
H^\hdot_F(X) @>>> H^\hdot_{DR}(X) @>>> H^\hdot(X,\calo_X) @>>> \dots
\end{CD}
$$
associated to the map $H^\hdot_{DR}(X) \to H^\hdot(X,\calo_X)$. If
$X$ is admissible, then the group $H^2_F(X)$ coincides with the
kernel of the natural map $H^2_{DR}(X) \to H^2(X,\calo_X)$.

\subsection{Definitions.}
The prototype for quantization is the quantization of a formal
polydisc. Let $\A$ be the power series algebra
$$
\A \cong k[[x_1,\ldots,x_d,y_1,\ldots,y_d]]
$$ 
on $2d$ variables $x_1,\ldots,x_d,y_1,\ldots,y_d$. Roughly speaking,
quantizing $\A$ consist of passing to the so-called {\em Weyl
algebra}.

\begin{defn}
The {\em formal Weyl algebra} (of fixed dimension $2d$) is the
complete topological associative algebra
$$
D=k[[x_1,\ldots,x_d,y_1,\ldots,y_d,h]]
$$
topologically generated by elements
$x_1,\ldots,x_d,y_1,\ldots,y_d,h$ subject to relations
\begin{align*}
[x_i,x_j] &= [y_i,y_j]=[x_i,h]=[y_j,h]=0,\\
[x_i,y_j] &=\delta_{ij}h
\end{align*}
for all $0 < i,j \leq d$.
\end{defn}

The formal Weyl algebra $D$ is a flat algebra over the power series
algebra $k[[h]]$. The subspace $hD \subset D$ is a two-sided ideal,
and the quotient $D/hD$ is isomorphic to the power series algebra
$\A$.

\medskip

The general definition of quantizations is as follows. Let $X$ be an
$S$-manifold with structure morphism $\pi:X \to S$, and denote by
$\pi^{-1}\calo_S$ the sheaf-theoretic pullback of the structure
sheaf $\calo_S$.

\begin{defn}\label{qua.defn}
A {\em quantization} $\D$ of an $S$-manifold $X$ is a sheaf of
associative flat $\pi^{-1}\calo_S[[h]]$-algebras on $X$ complete in
the $h$-adic topology and equiped with an isomorphism $\D/h\D \cong
\calo_X$.
\end{defn}

We note that quantizations are compatible with base change. Namely,
given an $S$-manifold $X$ with quantization $\D$ and a map $f:S' \to
S$, we obtain a quantization $f^*\D$ of the $S'$-manifold $X
\times_S S'$ by setting
$$
f^*\D = f^\hdot \D \whotimes_{f^\hdot\calo_S} \calo_{S'}.
$$
A particular case of this construction allows one to define jet
bundles for quantizations. Assume given a quantization $\D$ of an
$S$-manifold $X$. Consider the product $X \times_S X$ with the
projections $p_1,p_2:X \times_S X \to X$. The second projection
$p_2$ turns $X \times_S X$ into an $X$-manifold. Let $\D' = p_1^*\D$
be the quantization of the $X$-manifold $X \times_S X$ obtained by
pullback with respect to the projection $p_1$. Then $\D'$ is a sheaf
of $p_2^\hdot\calo_X$-algebras on $X \times_S X$, and we have
$\D'/h\D' \cong \calo_{X \times_S X}$. The ideal $\J_\Delta \subset
\calo_{X \times_S X}$ of the diagonal $X \cong \Delta \subset X
\times_S X$ lifts to a well-defined two-sided ideal $h\D' +
\J_\Delta \subset \D'$. The completion $J^\infty\D$ of the sheaf of
algebras $\D'$ with respect to the sheaf of ideals $h\D' +
\J_\Delta$ is supported on the diagonal, and it is naturally a sheaf
of $\calo_X$-algebras. Moreover, it is easy to see that $J^\infty\D$
is a pro-vector bundle on $X$. The fiber $J^\infty\D_x$ of the
bundle $J^\infty\D$ at a closed point $x \in X$ is canonically
isomorphic to the completion $\wh{\D_x}$ of the stalk $\D_x$ of the
sheaf $\D$ at the point $x \in X$ with respect to the topology
generated by the ideal $h\D_x + \m_x$, where $\m_x \subset \calo_x
\cong \D_x/h\D_x$ is the maximal ideal in the local ring $\calo_x$
of germs of functions on $X$ near $x \in X$.

\begin{defn}
The bundle $J^\infty\D$ is called the {\em jet bundle} of the
quantization $\D$.
\end{defn}

Quantizations are usually studied in connection with Poisson
geometry (see e.g. \cite{Kn}); we briefly recall this connection.
Given a quantization $\D$ on an $S$-manifold $X$, one considers the
commutator in the non-commutative algebra $\D$ and defines a
skew-symmetric bracket operation
$$
\{-,-\}: \calo_X \otimes \calo_X \to \calo_X
$$
by $\{a,b\} = \frac{1}{h}\wt{a}\wt{b}-\wt{b}\wt{a} \mod h^2$ for any
two local sections $a$, $b$ of the sheaf $\calo_X$ lifted to local
sections $\wt{a}$, $\wt{b}$ of the sheaf $\D$. One checks easily
that this bracket is well-defined and satisfies the axioms of a
Poisson bracket, namely,
$$
\{a,(b\cdot c)\} = \{a,b\}\cdot c + \{a,c\}\cdot b, \qquad
\{a,\{b, c\}\} + \{b,\{c,a\}\} + \{c,\{a,b\}\} = 0.
$$
By definition, this means that $X$ becomes a so-called {\em Poisson
scheme} over $S$, and one says that $\D$ is a quantization of the
Poisson scheme $X$.

Since a Poisson bracket $\{-,-\}$ is a derivation with respect to
both parameters, it is given by
$$
\{a,b\} = da \wedge db \cntrct \Theta
$$
for some bivector field $\Theta \in \Lambda^2\T_{X/S}$, where
$\T_{X/S}$ is the relative tangent bundle of $X$ over $X$. The
bivector field $\Theta$ defines an $\calo_X$-valued pairing on the
relative cotangent bundle $\Omega^1(X/S)$.

In this paper, we will only be interested in Poisson brackets such
that the associated pairing on $\Omega^1(X/S)$ is
non-degenerate. Since the pairing is skew-symmetric, this in
particular means that the dimension $\dim X/S$ must be
even. Applying the non-degenerate pairing $\Theta$, one identifies
$\T_{X/S}$ and $\Omega^1(X/S)$, so that $\Theta$ induces a $2$-form
$\Omega \in \Omega^2(X/S)$. Conversely, given a non-degenerate
$2$-form $\Omega \in \Omega^2(X/S)$, one applies it to identify
$\Omega^1(X/S)$ with $\T_{X/S}$ and obtains a non-degenerate
bivector field $\Theta \in \Lambda^2\T_{X/S}$. It is well-known that
$\Theta$ defines a Poisson bracket if and only if $\Omega$ is a
closed form. Thus giving a Poisson structure on an $S$-manifold $X$
with non-degenerate pairing $\Theta$ is the same as giving a
symplectic form $\Omega \in \Omega^2(X/S)$. Given a symplectic
$S$-manifold $X$ (of some even dimension $2d$), by a quantization of
$X$ we will understand a quantization of the $S$-manifold $X$ such
that the associated Poisson bracket on $X$ coincides with the
bracket induced by the symplectic form.

The definition of quantizations generalizes verbatim to the case of
formal schemes; in particular, it applies to the formal polydisc
$\Spf \A$ over a field $k$. Set-theoretically, $\Spf \A$ is a point,
so that a quantization $\D$ of $\Spf\A$ is an algebra over $k$. One
example of such a quantization is the formal Weyl algebra $D$. Our
approach to quantizations is based on the following standard fact
(essentially, a version of the Darboux Theorem).

\begin{lemma}\label{drboux}
Let $\D$ be any quantization of the formal polydisc $\Spf \A$ over a
field $k$ of characteristic $0$ such that the associated Poisson
pairing $\Theta$ is non-degenerate. Then $\D$ is isomorphic to the
formal Weyl algebra $D$.\endproof
\end{lemma}

In particular, for any quantization $\D$ of a smooth symplectic
manifold $X$ over a field of characteristic $0$, the completion
$\wh{\D}_x$ of the stalk $\D_x$ at some closed point $x \in X$ is a
quantization of the formal neighborhood of $x$ in $X$, which is
isomorphic to the formal polydisc $\Spf \A$ over the residue field
$k = k_x$ of the point $x \in X$. By Lemma~\ref{drboux}, there
exists a (non-canonical) isomorphism $\wh{\D}_x \cong D$.

\begin{remark}
Assume that both the base $S = \Spec O_S$ and a smooth $S$-manifold
$X = \Spec O_X$ are affine. Then for every quantization $\D$ of the
manifold $X$, the algebra $D_X=H^0(X,\D)$ of global sections of the
sheaf $\D$ is a flat $O_S[[h]]$-algebra, complete in the $h$-adic
topology and equiped with an isomorphism $D_X/hD_X \cong O_X =
H^0(X,\calo_X)$ (the natural map $D_X \to O_X$ is surjective because
the sheaf $\calo_X$ has no cohomology). Conversely, there exist a
non-commutative localization procedure called {\em Ore localization}
which applies, in particular, to any one-parameter deformation of
the algebra $O_X$ (see e.g. \cite[\S 2.1]{Kapranov}) and gives a
one-parameter deformation of the sheaf of algebras $\calo_X$ on
$X$. The constructions are mutually inverse. Thus in the affine
case, giving a quantization of $X=\Spec O_X$ is equivalent to giving
a one-parameter deformation of the $O_S$-algebra $O_X$ (more
precisely, a flat $O_S[[h]]$-algebra $D_X$ complete in the $h$-adic
topology and equiped with an isomorphism $D_X/hD_X \cong O_X$).
\end{remark}

\begin{remark}
In general, a quantization $\D$ of an $S$-manifold $X$ does not have
to be isomorphic to $\calo_X[[h]]$ even as a sheaf of
groups. However, this is true in the affine case over a point, $S =
\Spec k$, $k$ a field, $X = \Spec O_X$. Indeed, deformations of the
algebra $O_X$ in the class of associative $k$-algebras are
controlled by the so-called {\em Hochschild cohomology groups}
$HH^\hdot_i(O_X)$, which are computed by means of the Hochschild
cochain complex $C^\hdot_i(O_X)$,
$$ 
C^i_k(O_X) = \Hom_k(O_X^{\otimes i},O_X),
$$
where the tensor product is taken in the category of
$k$-modules. Inside $C^i_k(O_X)$, one distinguishes a subcomplex
$$
C^i_{diff}(O_X) \subset C^i_k(O_X)
$$
of cochains given by polydifferential operators. For any
quantization $\D$, the algebra $D_X$ is isomorphic to $k[[h]]$ as a
$k$-vector space; it is in the multiplication operation in $D_X$
that the non-triviality of the quantization is contained (this
multiplication is usually referred to as the {\em star-product}).
The complex $C^\hdot_{diff}(O_X)$ controls those deformations for
which the star-product is given by a series with bidifferential
operators as coefficients. However, an easy computation shows that
the complexes $C^\hdot_{diff}(O_X)$ and $C^\hdot_k(O_X)$ are
quasiisomorphic. Therefore any deformation $D_X$ in fact can be
given by a star-product whose Taylor coefficients in $h$ are
bidifferential operators. Since polydifferential operators are
local, -- that is, induced by sheaf maps $\calo_X^{\otimes k} \to
\calo_X$, -- after localization we can represent $\D$ as the sheaf
$\calo_X[[h]]$ with some non-trivial star-product multiplication.
\end{remark}

\subsection{Statements.}
We can now formulate our main result.

\begin{theorem}\label{main}
Let $X$ be an admissible $S$-manifold of dimension $2d$ equiped
with a closed non-degenerate relative form $\Omega \in
H^0(X,\Omega^2_{X/S})$. Denote by $[\Omega] \in H^2_{DR}(X/S)$ the
cohomology class of the symplectic form. Let $Q(X,\Omega)$ be the
set of isomorphism classes of quantizations of $X$ compatible with
the form $\Omega$.

Then there exists a natural injective map
$$
\Per:Q(X,\Omega) \hookrightarrow H^2_{DR}(X/S)[[h]],
$$ 
called the {\em non-commutative period map}. Moreover, for every
quantization $q \in Q(X,\Omega)$, the power series $f = \Per(q) \in
H^2_{DR}(X)[[h]]$ has constant term $[\Omega]$. Finally, any
splitting $P:H^2_{DR}(X/S) \to H^2_F(X/S)$ of the canonical
embedding $H^2_F(X/S) \to H^2_{DR}(X/S)$ induces an isomorphism
$$
P \circ \Per:Q(X,\Omega) \quad \overset{\sim}{\longrightarrow} \quad
P([\Omega]) + h H^2_F(X/S)[[h]] \subset H^2_F(X/S)[[h]]
$$
between $Q(X,\Omega)$ and the set of all power series in $h$ with
coefficients in $H^2_F(X/S)$ and constant term $P([\Omega])$.
\end{theorem}

In particular, quantizations always exist (provided the manifold in
question is admissible). Moreover, one can define a preferred
quantization:

\begin{defn}\label{can.defn}
A quantization $\D \in Q(X,\Omega)$ of an admissible symplectic
$S$-manifold is called {\em canonical} if its period $\Per(\D) \in
H^2_{DR}(X/S)[[h]]$ is the constant power series $[\Omega]$.
\end{defn}

The period map itself is completely canonical. However, the
pa\-ra\-me\-tri\-zation of quantizations by formal power series with
coefficients in $H^2_F(X/S)$ does depend on the splitting
$F:H^2_{DR}(X/S) \to H^2_F(X/S)$. Sometimes there is a canonical
choice of this splitting -- for instance, when $X$ is projective
over $\C$, such a splitting is provided by Hodge theory. The
canonical quantization enjoys several nice properties, but we should
warn the reader that it {\em does not} have to exist, even for an
admissible manifold -- unless $H^i(X,\calo_X)=0$ for $i=1,2$, so
that the period map is surjective.


In the basic case when $S = \Spec k$ is a point, Theorem~\ref{main}
is completely parallel to what one has for symplectic deformations
-- in other words, for commutative deformations of the pair $\langle
X,\Omega \rangle$. The commutative version of the period map was
introduced in \cite{KV}, and it is very simple: it sends a
deformation to the associated family of cohomology classes
$[\Omega]_h \in H^2_{DR}(X)$ of the corresponding symplectic
forms. This motivates our terminology. Unfortunately, we do not have
a similar interpretation of the non-commutative period map.

Our definition of the period map is also quite simple, in fact, it
takes one paragraph --- and the full proof of Theorem~\ref{main}
takes only two pages. Both are contained in Section~\ref{per}. But
both the definition and the proof require some preliminary
machinery. All the facts we need are essentially standard, but there
are no suitable references in the literature. Thus we have to devote
Section~\ref{hc} and Section~\ref{fg} to these preliminaries. So as
not to overwhelm the reader with technicalities, some proofs are
postponed till Section~\ref{na} (which only depends on
Section~\ref{hc}). Section~\ref{ex} contains some extensions of our
results to other frameworks. In particular, we consider the
equivariant version of Theorem~\ref{main}. We also clarify the
relation between quantizations and the universal symplectic
deformation constructed in \cite{KV} by showing that both can be
incorporated into a single multi-parameter partially non-commutative
deformation. It is here that the general relative setting of
Theorem~\ref{main} plays a crucial role. Finally, Section~\ref{dis}
is taken up with some concluding remarks --- we try to place our
results in the general context and compare them with existing
alternative approaches to deformation quantization.

\section{Preliminaries on Harish-Chandra torsors.}\label{hc}

\subsection{Harish-Chandra pairs.}
The following definition was first introduced most probably by
A. Beilinson and J. Bernstein, \cite{BB}.

\begin{defn}\label{hc.defn} 
A {\em Harish-Chandra pair} $\GH$ over the field $k$ is a pair of a
  connected affine algebraic group $G$ over $k$, a Lie algebra $\h$
  over $k$ equiped with a $G$-action, and an embedding $\g \to \h$
  of the Lie algebra $\g$ of the group $G$ into the Lie algebra $\h$
  such that the adjoint action of $\g$ on $\h$ is the differential
  of the given $G$-action.

A {\em module} $V$ over a Harish-Chandra pair $\langle G, \h
\rangle$ is a representation $V$ of the Lie algebra $\h$ whose
restriction to $\g \subset \h$ is intergrated to an algebraic
representation of the group $G$.
\end{defn}

Just as when working with jet bundles, in applications it is
important to allow groups which are not finite-dimensional, or, more
precisely, to allow $G$ to be the projective limit of affine
algebraic groups. To extend Definition~\ref{hc.defn} to this case,
we make the following modifications. The Lie algebra $\g$ is a
topological vector space equiped with a ``compact'' topology --
namely, it is a projective limit of finite-dimensional vector
spaces. Note that topological vector spaces of this type form an
abelian category (the one dual to the category of usual vector
spaces). The Lie algebra $\h$ is also a projective limit of
finite-dimensional vector spaces and moreover, $\g \subset \h$ is of
finite codimension (in other words, $\g$ is closed in $\h$). The
group $G$ will always be an affine group scheme and a projective
limit of affine algebraic groups of finite type over $k$.

A module $V$ over a Harish-Chandra pair will also be a projective
limit of finite-dimensional vector spaces, and we will assume that
both $G$ and $\h$ act in a way compatible with this topology. In the
case of the group scheme $G = \Spec A$, this means that the
$G$-action on $V$ is given by a coaction $V^* \to V^* \otimes_k A$
of the Hopf algebra $A$ on the (discrete, although
infinite-dimensional) vector space $V^*$ topologically dual to
$V$. Modules defined in this way also form an abelian category. This
category comes equiped with a symmetric tensor product (defined in
the obvious way). The unit object for this product is the
one-dimensional trivial representation $k$.

As usual both for groups and for Lie algebras, given a $\GH$-module
$V$, by the {\em cohomology groups} $H^\hdot(\GH,V)$ of the module
$V$ we will understand the $\Ext$-groups $\Ext^\hdot(k,V)$ (taken in
the category of topological Harish-Chandra modules).

\begin{remark}
There is a more general notion of a Harish-Chandra pair (see
\cite{BFM}), where the Lie algebra $\h$ is allowed to be a so-called
{\em Tate topological vector space}. In this paper, we do not need
it.
\end{remark}

\subsection{Torsors.}
Let $X$ be an $S$-manifold. To keep things precise, we will say that
given a group scheme $G$, by a $G$-torsor over $X$ we will
understand a scheme $Y$ faithfully flat over $X$ and equiped with
an action map $G \times Y \to Y$ which commutes with the projection
to $X$ and induces an isomorphism $G \times Y \to Y \times_X Y$. (In
our applications, all torsors will be locally trivial in Zariski
topology.) Assume given a Harish-Chandra pair $\GH$. For any
$G$-torsor $\M$ over $X$ we have the Lie algebra bundles $\g_\M$ and
$\h_\M$ on $X$ associated to the $G$-modules $\g$ and $\h$. The map
$\g \to \h$ induces a map $\g_\M \to \h_M$. Moreover, since we work
in characteristic $0$, the scheme $G$ is smooth, so that the
faithfully flat projection $\rho:\M \to X$ is also smooth. Therefore
we have a $G$-equivariant short exact sequence
$$
\begin{CD}
0 @>>> \T_{\M/X} @>>> \T_{\M/S} @>>> \rho^*\T_{X/S} @>>> 0
\end{CD}
$$
of relative tangent bundles, which by descent gives the so-called
{\em Atiyah extension}
\begin{equation}\label{ati}
\begin{CD}
0 @>>> \g_M @>{\iota_\M}>> \E_\M @>>> \T_{X/S} @>>> 0
\end{CD}
\end{equation}
of bundles on $X$. Recall that a $G$-invariant
connection\footnote{As noted in Subsection~\ref{dfn}, connections
are understood relatively over $S$.} on the principal $G$-bundle
$\M$ is by definition given by a bundle map $\theta_\M:\E_\M \to
\g_M$ which splits the extension \eqref{ati} -- in other words, the
composition $\theta_\M \circ \iota_\M:\g_\M \to \g_\M$ is the
idenitity map. Equivalently, one can specify the corresponding
$G$-invariant $\g$-valued $1$-form $\rho^*\theta_\M \in
H^0(\M,\Omega^1_\M \otimes \g)$, where $\rho:\M \to X$ is the
projection.  A connection $\theta_\M$ is {\em flat} if the
corresponding one-form $\Omega = \rho^*\theta_\M$ satisfies the
Maurer-Cartan equation $2d\Omega + \Omega \wedge \Omega =
0$. Generalizing this, by an {\em $\h$-valued connection} on $\M$ we
will understand a bundle map $\theta_\M:\E_M \to \h_\M$ such that
the composition $\theta_\M \circ \iota_\M:\g_\M \to \h_\M$ is the
given embedding. Again, an $\h$-valued connection $\theta_\M$ is
called {\em flat} if the corresponding $\h$-valued $1$-form $\Omega
= \rho^*\theta_\M$ satisfies $2d\Omega + \Omega \wedge \Omega = 0$.

\begin{defn} By a {\em Harish-Chandra $\GH$-torsor} $\M$ over the
$S$-manifold $X$ we will understand a pair $\langle
\M,\theta_\M\rangle$ of a $G$-torsor $\M$ over $X$ and a flat
$\h$-valued connection $\theta_\M:\E_M \to \h_M$ on $\M$.
\end{defn}

The notion of a Harish-Chandra torsor has the usual
functorialities. In particular, if we have a map of Harish-Chandra
pairs $f:\GH \to \GHh$ and a $\GH$-torsor $\M$, then we canonically
obtain the induced $\GHh$-torsor $\M_1 = f_*\M = \M \times^G
G_1$. For a tautological Harish-Chandra pair $\langle G,\g \rangle$,
a $\langle G,\g \rangle$-torsor over $X$ is the same as a principal
$G$-bundle equiped with a $G$-invariant flat connection.

The set of isomorphism classes of all $\GH$-torsors over an
$S$-manifold $X$ will be denoted by $H^1(X,\GH)$. The torsors
themselves form a category. This category is a groupoid, which we
will denote by $\hh^1(X,\GH)$.

An important special class of Harish-Chandra torsors is the
following one.

\begin{defn}\label{trans}
A $\GH$-torsor $\M$ over an $S$-manifold $X$ is called {\em
transitive} if the connection map $\theta_\M:\E_\M \to \h_M$ is an
isomorphism.
\end{defn}

Given a transitive $\GH$-torsor $\M$, one can invert the connection
map and obtain a bundle map $\theta_\M^{-1}:\h_\M \to \E_\M$, or,
equivalently, a $G$-equivariant map $\h \otimes \calo_\M \to
\T_{\M/S}$ of vector bundles on $\M$. The latter is in turn
equivalent to giving a $G$-equivariant map $\h \to
H^0(\M,\T_{\M/S})$, and it is easy to check that this is a Lie
algebra map if and only if the connection $\theta_\M$ is flat.  Thus
for a transitive $\GH$-torsor $\M$, the whole Harish-Chandra pair
$\GH$ acts on the scheme $\M$ -- that is, the Lie algebra $\h$ acts
by derivations of the structure sheaf, and the restriction of this
action to $\g \subset \h$ is the differential of the $G$-action. The
$\h$-action is transitive, which explains our
terminology. Conversely, a $G$-torsor $\M$ equiped with a
compatible transitive action of $\h$ gives rise to a transitive
$\GH$-torsor in the sense of Definition~\ref{trans}.

\subsection{Localization.}
Assume given an $S$-manifold $X$, a Harish-Chandra pair $\GH$, and a
$\GH$-torsor $\langle\M,\theta_\M\rangle$ over $X/S$. Let $V$ be a
finite-dimensional $\GH$-module. Then we have a map
$$
f:\GH \to \langle GL(V),\gl(V) \rangle
$$ 
and the induced torsor $f_*\M$. If $\V$ is the vector bundle on $X$
associated to the $G$-module $V$, then $f_*\M$ coincides with the
principal $GL(V)$-bundle of frames in $\V$. By construction it
carries a flat connection. Thus $\V$ also carries a canonical flat
connection $\nabla$. Explicitly, let $\xi \in \Gamma(U,\E_\M)$ be a
local section of the Atiyah sheaf $\E_\M$, and let $a \in
\Gamma(U,\V)$ be a local section of the bundle $\V$. Then by
construction both the Atiyah sheaf and the Lie algebra bundle
$\h_\M$ act on sections of the bundle $\V$, and the expression
\begin{equation}\label{loc}
\nabla_\xi(a) = \xi \cdot a - \theta_\M(\xi) \cdot a \in \Gamma(U,\V)
\end{equation}
only depends on the image of $\xi$ in the tangent sheaf $\T_X$. Thus
it defines a connection on $\V$, which is exactly $\nabla$.

When the module $V$ is only a projective limit of finite-dimensional
vector spaces, the group $GL(V)$ is not well-defined. However, we
can still define a flat connection on the associated bundle $\V$ by
directly applying \eqref{loc}. Associated bundle in this case lies
in the completed category of pro-coherent sheaves -- just as the jet
bundles considered in Subsection~\ref{dfn}.

To sum up, given the torsor $\M$, for any module $\GH$-module $V$ we
obtain a flat bundle $\V$ on the $S$-manifold $X$. In other words,
the torsor $\M$ defines a functor from the category of $\GH$-modules
to the category of flat vector bundles on $X/S$. We will call this
the {\em localization functor} associated to $\M$, and we will
denote it by
$$
\V = \loc(\M,V).
$$ 
The functor of localization with respect to $\M$ is obviously
exact. In particular, it extends to derived categories and induces a
canonical localization map
$$
\loc(\M,-):H^\hdot(\GH,V) \to H^\hdot_{DR}(X,\V).
$$
Moreover, localization is a tensor functor.

\begin{remark}
Most probably, the converse is also true: modulo the appropriate
finiteness conditions, every tensor functor from the category of
$\GH$-modules to the category of flat bundles on $X$ comes from a
$\GH$-torsor $\M$ on $X$. Equivalent functors give isomorphic
torsors. We do not develop this Tannakian-type formalism here to
save space.
\end{remark}

Localization can be also be described in a different
language. Recall that the standard descent procedure induces an
equivalence between the category of vector bundles on $X$ equiped
with a flat connection and the category of $G$-equivariant vector
bundle on $\M$ equiped with a flat connection which is compatible
with the $G$-action. Compatibility here means that for every vector
$\xi \in G$, the covariant derivative $\nabla_\xi$ with respect to
the corresponding vector field on $\M$ coincides with the action of
$\xi$ coming from the $G$-equivariant structure (it is well-known
that this definition does not give the correct equivariant version
of the derived category of flat vector bundles; however, if we stick
to the abelian categories, the descent works just fine). Using this
equivalence, one does the localization procedure in two
steps. First, one considers the constant vector bundle $V \otimes
\calo_\M$ on $\M$ with the trivial flat connection, and equips it
with the product $G$-action. The connection and the $G$-action are
not compatible. Then one corrects the connection on $V \otimes
\calo_\M$ by \eqref{loc} -- for this one needs the $\h$-valued
connection on $\M$ and the $\h$-action on $V$. After that, the
localization $\loc(\M,V)$ is obtained by descent.

The descent procedure is of course quite general, it is by no means
limited to vector bundles of type $V \otimes \calo_\M$. We note that
descent works especially well when the $\GH$-torsor $\M$ is
transitive in the sense of Definition~\ref{trans}. In this case, we
have an $\h$-action on $\M$ which trivializes the tangent bundle
$\T_\M$, $\T_\M \cong \h \otimes \calo_\M$ (the Lie algebra
structure on $\h \otimes \calo_\M$ is not $\calo_\M$-linear, it is
skew-linear with respect to the $\h$-action on $\calo_\M$). It
follows immediately that giving a compatible flat connection
$\nabla$ on a $G$-equivairant vector bundle $\E$ on $\M$ is
equivalent to extending the $G$-action on $\E$ to a compatible
$\h$-action. Thus we have the following.

\begin{lemma}\label{desce}
Let $\M$ be a transitive $\GH$-torsor over an $S$-manifold $X$. Then
the category of vector (pro)bundles on $X$ equiped with a flat
connection is equivalent to the category of $\GH$-equivariant vector
(pro)bundles on $\M$.\endproof
\end{lemma}

Here a $\GH$-equivariant vector bundle $\E$ is a $G$-equivariant
vector bundle equiped with an action of the Lie algebra $\h$ which
is compatible with the $\h$-action on $\M$ and gives the
differential of the $G$-action after restriction to $\g \subset \h$.
For a constant $\GH$-module $V$, one simply takes the product
$\GH$-action on $V \otimes \calo_\M$; descent by Lemma~\ref{desce}
gives $\loc(\M,V)$.

\subsection{Harish-Chandra extensions.}\label{nas} 
In the body of the paper, we will need to study the behavior of
Harish-Chandra torsors under extensions. More precisely, we need
what is usually referred to as the long exact sequence in the
non-abelian cohomology. So as not to interrupt the exposition too
much, we give all the statements here, and we postpone the proofs
till Section~\ref{na}.

Let $\GH$ be a Harish-Chandra pair, and let $V$ be a $\GH$-module.
Consider $V$ as an (additive) algebraic group. By an {\em extension}
\begin{equation}\label{extt}
\begin{CD}
1 @>>> V @>{\rho}>> \GHh @>{\pi}>> \GH @>>> 1
\end{CD}
\end{equation}
of the pair $\GH$ by the module $V$ we will understand a
Harish-Chandra pair $\GHh$ equiped with a map $f:\GHh \to \GH$ such
that $\Ker f = \vv \subset \GHh$ is the tautological Harish-Chandra
pair associated to $V$, and the adjoint action of $\GH$ on $V$ comes
from the given module structure. In other words, we have an
extension of groups compatible with the extension of the Lie
algebras.

Given an $S$-manifold and a $\GH$-torsor $\M$ over $X/S$, we denote
by 
$$
H^1_\M(X,\GHh)
$$
the set of isomorphism classes of $\GHh$-torsors $\M_1$ on $X/S$
equiped with an isomorphism $\pi_*\M_1 \cong \M$. We will call
torsors of this type {\em liftings} of the torsor $\M$ to the
Harish-Chandra pair $\GHh$ (if $\GHh$ were to be a subobject $\GHh
\subset \GH$, the common term would be ``restriction'').

Let $\V = \loc(\M,V)$ be the localization of the $\GH$-module
$V$ 
with respect to the torsor $\M$. The basic statement we need is the
following one.

\begin{prop}\label{nab}\mbox{}
Let $\M$ be a transitive $\GH$-torsor over an $S$-manifold $X$.
\begin{enumerate}
\item There exists a canonical cohomology class $c \in H^2(\GH,V)$
with the following property: the set $H^1_\M(X,\GHh)$ is non-empty
if and only if the localization $\loc(\M,c) \in H^2_{DR}(X/S,\V)$ is
trivial.
\item If the class $\loc(\M,c)$ is indeed trivial, then the set
$H^1_\M(X,\GHh)$ is naturally a torsor over the de Rham cohomology
group $H^1_{DR}(X/S,\V)$.
\end{enumerate}
\end{prop}

We will also need a more involved statement, a certain compatibility
result vaguely reminiscent of the octahedron axiom in homological
algebra. Consider a Harish-Chandra pair $\GH$, and let
\begin{equation}\label{ab}
\begin{CD}
0 @>>> U @>{a}>> V @>{b}>> W @>>> 0
\end{CD}
\end{equation}
be a short exact sequence of $\GH$-modules. Assume given an extension
$$
\begin{CD}
1 @>>> V @>>> \GHh @>{\pi}>> \GH @>>> 1
\end{CD}
$$
of the Harish-Chandra pair $\GH$ by the module $V$, and denote its
cohomology class by $c \in H^2(\GH,V)$. Let $\GHt = \GHh/U$ be the
associated extension of $\GH$ by $W$. By definition, $\GHh$ is an
extension of $\GHt$ by the module $U$. Denote its cohomology class
by $c_0 \in H^2(\GHt,U)$.

Assume given a $\GH$-torsor $\M$ over $X/S$, and let $\U$, $\V$ and
$\W$ be the localizations of the $\GH$-modules $U$, $V$ and $W$. We
have a long exact sequence of de Rham cohomology groups
\begin{equation}\label{lng}
\begin{CD}
H^2_{DR}(X/S,\U) @>{a}>> H^2_{DR}(X/S,\V) @>{b}>> H^2_{DR}(X/S,\W)
@>>>
\end{CD}
\end{equation}
Assume that we are in the following situation: the localization
$$
\loc(\M,c) \in H^2_{DR}(X/S,\V)
$$ 
is {\em not} trivial, but its restriction
$$
b(\loc(\M,c)) \in H^2_{DR}(X/S,\W)
$$ 
{\em is} trivial. Then the $\GH$-torsor $\M$ does not lift to a
$\GHh$-torsor over $X/S$, but it does lift to a
$\GHt$-torsor. Moreover, for every such lifting $\M_1 \in
H^1_\M(X/S,\GHt)$, we obtain a lifting of the obstruction cohomology
class $\loc(\M,c) \in H^2_{DR}(X/S,\V)$ to a cohomology class in
$H^2_{DR}(X/S,\U)$, namely, the class $\loc(\M_1,c_0) \in
H^2_{DR}(X/S,\U)$. By Proposition~\ref{nab}, we know that the set
$H^1_\M(X/S,\GHt)$ is a torsor over the group $H^1_{DR}(X/S,\W)$. On
the other hand, by the exact sequence \eqref{lng} the group
$H^1_{DR}(X/S,\W)$ acts on the set $H^2_{DR}(X/S,\U)$.

\begin{lemma}\label{oct}
The map $H^1_\M(X/S,\GHt) \to H^2_{DR}(X/S,\U)$ given by
$$
\M_1 \mapsto \loc(\M_1,c_0)
$$
is compatible with the $H^1_{DR}(X/S,\W)$-action on both sides.
\end{lemma}

The reader will find the proofs of Proposition~\ref{nab} and
Lemma~\ref{oct} in Section~\ref{na}. We note that the condition of
transitivity is needed only in the proof of the ``if'' part of
Proposition~\ref{nab}~\thetag{i}. The Proposition is also true for
general $\GH$-torsors, but the proof is slightly longer; since in
our applications all torsors will be transitive, we impose this
assumption right from the start to save space.

\section{Quantization via formal geometry.}\label{fg}

\subsection{The bundle of coordinate systems.}\label{fgs}
Formal geometry is a technique of dealing with various questions in
differential geometry by solving them first in the universal
context, -- that is, over a formal polydisc, -- and equivariantly
with respect to the Lie algebra of vector fields on the polydisc. It
dates back at least to the papers \cite{GK} by I. Gelfand and
D. Kazhdan and/or \cite{bt} by R. Bott. However, there are no
convenient general references. We have learned what we know of this
technique at B. Feigin's Moscow seminar. Since we do need to use it,
-- and in the relative setting, to make things worse, -- we give
here a self-contained outline of the basic setup.

\bigskip

Fix a dimension $n$, at this point not necessarily even. Consider
the formal power series algebra $\A=k[[x_1,\ldots,x_n]]$. Denote by
$\WW$ the Lie algebra of derivations of the algebra $\A$ -- in other
words, the Lie algebra of vector fields on the formal polydisc $\Spf
\A$. Consider the subalgebra $\WW_0 \subset \WW$ of vector fields
vanishing at the closed point (equivalently, derivations preserving
the maximal ideal in $\A$). Then the Lie algebra $\WW_0$ is
naturally the Lie algebra of a proalgebraic group $\Aut \A$ of
automorphisms of the local $k$-algebra $\A$. In the language of
Section~\ref{hc}, we have a Harish-Chandra pair $\WA$.

Let $X$ be an $S$-manifold of dimension $n$ with projection map
$\pi:X \to S$. For any scheme $T$, giving a map $\eta: T \to X$ is
equivalent to giving a map $p(\eta) = \pi \circ \eta:T \to S$ and a
section $\sigma(\eta):T \to T \times_S X$ of the canonical
projection $T \times_S X \to T$. Formal germs of functions on $T
\times_S X$ near the closed subscheme $\sigma(\eta)(T) \subset T
\times_S X$ form a sheaf of (topological) $\calo_T$-algebras on $T$
which we denote by $\wh{\calo}_{X,\eta}$. If the scheme $T$ is
affine, the sheaf $\wh{\calo}_{X,\eta}$ is non-canonically
isomorphic to the completed tensor product $\calo_T \whotimes
\A$. Let $\M(T)$ be the set of all pairs
$$
\left\langle \eta:T \to X, \quad \phi:\wh{\calo}_{X,\eta} \cong
\calo_T \whotimes \A\right\rangle,
$$
where $\phi$ is an isomorphism of sheaves of topological
$\calo_T$-algebras. Geometrically, such a pair corresponds to a
commutative diagram
\begin{equation}\label{coord}
\begin{CD}
\Spf \A \times T         @>{\phi}>>     X     \\
       @VVV              @VV{\pi}V \\
         T              @>{p(\eta)}>>        S
\end{CD}
\end{equation}
which induces an identification between the formal neighborhood of
$\sigma(\eta)(T)$ in $T \times_S X$ and the product $\Spf \A \times
T$ -- loosely speaking, a family of formal coordinate systems on
$X/S$ parametrized by $T$. Setting $T \mapsto \M(T)$ defines a
functor from the category of affine schemes to the category of
sets. We leave it to the reader to check that this functor is
represented by a (non-Noetherian) scheme $\M_{coord}$, smooth and
affine over $X$. In fact, $\M_{coord}$ is the projective limit of a
family of $S$-manifolds, and it is a torsor over the group $\Aut\A$
with respect to the natural action. Moreover, the torsor
$\M_{coord}$ carries a structure of a transitive Harish-Chandra
torsor over $\WA$. Indeed, the Lie algebra $\WW$ also acts on $\A$,
hence on $\M_{coord}$, and the action map descends to a map
$$
a:\W_\M \to \E_\M,
$$ 
where $\W_\M$ is the vector bundle on $X$ associated to $\WW$, and
$\E_\M$ is the Atiyah bundle of the torsor $\M_{coord}$. It is
elementary to check that the map $a$ is in fact an isomorphism. To
define a $\WW$-valued flat connection $\theta_\M:\E_\M \to \W_\M$ on
$\M_{coord}$, it suffices to take the inverse isomorphism $\theta_\M
= a^{-1}$. Since $a$ is obtained from a Lie algebra map $\WW \to
\T_\M$, the connection $\theta_\M$ is flat.

\begin{defn}\label{coord.defn}
The $\WA$-torsor $\langle \M_{coord},\theta_\M = a^{-1}\rangle$ over
$X$ is called the {\em bundle of formal coordinate systems} on the
$S$-manifold $X$.
\end{defn}

The bundle of formal coordinate systems is the main object of formal
geometry. It is completely canonical, and it allows one to do the
following two things:
\begin{enumerate}
\item Obtain various canonical sheaves on $X$, such as sheaves of
sections of different symmetric and tensor powers of the tangent
bundle $\T(X)$, as sheaves of flat sections of localizations of
appropriate representations of the Harish-Chandra pair $\WA$.
\item Describe various differential-geometric structures on $X$ as
reductions of the torsor $\M_{coord}$ to different subgroups in
$\WA$.
\end{enumerate}
Usual applications revolve around \thetag{i}. More precisely, the
construction one uses is the following one. The simplest module over
the Harish-Chandra pair $\WA$ is the algebra $\A$ itself. It is easy
to check that its localization with respect to the $\WA$-torsor
$\M_{coord}$ coincides with the jet bundle $J^\infty\calo_X$:
$$
\loc(\M_{coord},\A) \cong J^\infty\calo_X.
$$
The sheaf of its flat sections is the structure sheaf $\calo_X$ of
the variety $X$. Analogously, one can take the $\WA$-module $\WW$ of
vector fields on $\A$, or the module $\Omega^p\A$ of $p$-forms on
$\A$ for some $p \leq n$, or, more generally, the $\WA$-module $\Xi$
of tensors of some type on $\A$. Then its localization is
$J^\infty\T$, the jet bundle of the tangent sheaf $\T_X$,
resp. $J^\infty\Omega^p_X$, resp. the jet bundle of the sheaf of
tensors on $X$ of the same type as $\Xi$. As usual, one recovers the
sheaf from its jet bundle by taking flat sections.

One can use this construction, for instance, to obtain
characteristic classes of the variety $X$ starting from cohomology
classes of the $\WA$-module $\Xi$. In the present paper, we leave
this subject completely alone. Our applications of formal geometry
are related to \thetag{ii}.

The following is the motivating example. Assume that the dimension
$n = 2d$ is even, and equip the formal polydisc $\A =
k[[x_1,\ldots,x_d,y_1,\ldots,y_d]]$ with the symplectic form $\omega
= \sum dx_i \wedge dy_i$. Denote by $\HH \subset \WW$ the Lie
subalgebra of Hamiltonian vector fields -- in other words, the
vector fields that preserve the symplectic form. As before, the
subalgebra $\WW_0 \cap \HH \subset \HH$ is naturally integrated to a
pro-algebraic group $\Sympl \A$, and we have a Harish-Chandra pair
$\WSA$.

\begin{lemma}\label{drb}
Let $X$ be an $S$-manifold of dimension $n=2d$.  There is a
one-to-one correspondence between symplectic structures on $X/S$ and
reductions of the $\WA$-torsor $\M_{coord}$ to $\WSA \subset \WA$.
\end{lemma}

\proof{} Given a symplectic form, one takes the subvariety $\M_s
\subset \M_{coord}$ of formal coordinate systems $\phi:\Spf\A \to X$
compatible with symplectic forms on both sides, and notices that by
the formal Darboux Theorem, $\M_s$ is a torsor over $\Sympl\A
\subset \Aut\A$.

Conversely, given such a reduction $\M_s \subset \M_{coord}$, one
recalls that by definition, the $\WA$-module $\Omega^2\A$ contains
an $H$-invariant vector $\omega$. By localization, $\omega$ gives a
flat section of the jet bundle $J^\infty\Omega^2_X$, thus a
symplectic form on $X/S$.
\endproof

Given a symplectic $S$-manifold, we will use the term {\em bundle of
symplectic formal coordinate systems} for the associated
$\WSA$-torsor $\M_s$ on $X$.

The cohomology class $[\Omega] \in H^2_{DR}(X)$ of the symplectic
form $\Omega$ also has a natural interpretation in terms of the
torsor $\M_s$. Namely, the de Rham complex of the polydisc gives a
resolution of the trivial $\WSA$-module $k$, and the standard
symplectic form on the polydisc defines a cohomology class $[\omega]
\in H^2(\WSA,k)$. The class $[\Omega] \in H^2_{DR}(X)$ is the
localization of this class $[\omega]$ with respect to the torsor
$\M_s$. The class $[\omega]$ corresponds to the central extension
$$
\begin{CD}
0 @>>> k @>>> \A @>>> H @>>> 0
\end{CD}
$$
of the Lie algebra $H$ of Hamiltonian vector fields on the polydisc
(or rather, to the corresponding extension of Harish-Chandra pairs).

\begin{remark} V. Drinfeld has explained to us that the torsor $\M_{coord}$
of formal coordinate systems on a manifold $X$ over an algebraically
closed field $k$ can be in fact characterized by a universal
property (see \cite{BD}, especially Remark 2.6.4 and Example
2.6.5). Firstly, for any Harish-Chandra pair $\GH$, there exists a
canonical formal manifold, carrying a transitive $\GH$-torsor (in
the language of \cite{BD}, our transtitive $\GH$-torsors correspond
to {\em $(\h,G)$-structures}). This manifold is $X = \Spf
k[[(\h/\g)^*]]$, the formal completion at $0$ of the quotient
$\h/\g$. Secondly, for any Harish-Chandra pair $\GH$ and a
transitive $\GH$-torsor $\M$ on a manifold $X$ of dimension $n$,
there exists a unique map $\GH \to \WA$ of Harish-Chandra pairs and
a unique compatible map $\tau:\M \to \M_{coord}$ from $\M$ into the
torsor $\M_{coord}$. To construct the map $\tau$, one fixes an
isomorphism between $\Spf k[[(\h/\g)^*]]$ and the standard formal
polydisc $\Spf \A$ and notices that a point $m \in \M$ lying over a
point $x \in X$ induces an isomorphism between $\Spf k[[(\h/\g)^*]]$
and the formal neighborhood of $x \in X$. In particular, the torsor
$\M_{coord}$ is the unique, up to a unique isomorphism, transitive
$\WA$-torsor over $X$.
\end{remark}

\subsection{Automorphisms of the formal Weyl algebra.}
We can now reformulate the quantization problem in the language of
formal geometry; additionally, we will recall some standard facts on
automorphisms of the formal Weyl algebra $D$.

The basic statement is very straightforward; essentially, it is a
quantized version of Lemma~\ref{drb}. Consider the Lie algebra $\Der
D$ of $k[[h]]$-linear derivations of the $k[[h]]$-algebra $D$. Since
the derivations in $\Der D$ are $k[[h]]$-linear, the algebra $\Der
D$ it preserves the ideal $hD \subset D$ and acts on the quotient
$\A =D/hD$. The action map $a:(\Der D) \to \WW$ factors through the
quotient $(\Der D)_0$ which is isomorphic to the Lie algebra $\HH
\subset \WW$ of Hamiltonian vector fields on the polydisc. The
subalgebra $(\Der D)^0 = a^{-1}(\WW^0)$ is naturally integrated to
the pro-algebraic group $\Aut D$, namely, the group of
$k[[h]]$-linear automorphisms of the Weyl algebra $D$ preserving the
two-sided ideal $\m_\A + hD \subset D$. Thus we have a natural
Harish-Chandra pair $\WDA$.

\begin{lemma}\label{qdrb}
Let $X$ be an $S$-manifold of dimension $n=2d$ equiped with a
symplectic form $\Omega$. Denote by $\M_s$ the bundle of symplectic
formal coordinate systems on $\langle X,\Omega \rangle$.

Then there exists a natural bijection between the set $Q(X,\Omega)$
of isomorphism classes of quantizations of the symplectic
$S$-manifold $X$, and the set $H^1_{\M_s}(X,\WDA)$ of the
isomorphism classes of liftings of the symplectic coordinate system
bundle $\M_s$ to a $\WDA$-torsor with respect to the canonical map
of Harish-Chandra pairs $\WDA \to \WSA$.
\end{lemma}

\proof{} To pass from a lifting $\M_q$ to a quantization, one takes
the localization $\loc(\M_q,D)$ of the $\WDA$-module $D$, and
considers the sheaf $\D$ of its flat sections. Since $D$ is a
$k[[h]]$-algebra, $\A = D/hD$, and both these facts are
$G$-equivariant, the sheaf $\D$ is a quantization of the symplectic
manifold $\langle X/S, \Omega \rangle$ in the sense of
Definition~\ref{qua.defn}.

Conversely, given a quantization $\D$ and an affine scheme $T$, one
follows Definition~\ref{coord.defn} and defines $\M_q(T)$ to be the
set of all pairs of a map $\eta:T \to X$ and an isomorphism
$$
\Phi:\calo_T \whotimes D \cong \wh{\D}_{X,\eta},
$$ 
where $\wh{\D}_{X,\eta}$ is the completion of the quantization
$p(\eta)^*\D$ of the $T$-manifold $T \times_S X$ obtained by
pullback with respect to the composition $p(\eta):T \to S$ of the
map $\eta:T \to X$ and the projection $X \to S$; the completion is
taken with respect to the ideal spanned by $h(p(\eta)^*\D)$ and the
ideal $\J_\eta \subset \calo_{T \times_S X}$ of the closed subscheme
$\sigma{\eta}(T) \subset T \times_S X$. We claim that the functor $T
\mapsto \M_q(T)$ is represented by a (non-Noetherian) scheme $\M_q$.

Indeed, consider the jet bundle $J^\infty\D$ of the quantization
$\D$. By definition, the completion $\wh{\D}_{X,\eta}$ coincides
with the pullback $\eta^*J^\infty\D$. Then every isomorphism
$\Phi:\calo_T \otimes D \to \wh{\D}_{X,\eta} \cong \eta^*J^\infty\D$
is defined by the $2n$ elements $\Phi(x_1),\dots,\Phi(x_n),
\Phi(y_n), \dots, \Phi(y_n) \in H^0(T,\eta^*J^\infty\D)$ which come
from the generators $x_1,\dots,x_n, y_1,\dots,y_n \in D$.
Conversely, such a set of $2n$ elements give a map if and only if
they satisfy the defining relations for $D$. Moreover, by the very
nature of these defining relations, every such map $\Phi$ induces a
symplectic map $\overline{\Phi}:\wh{\calo}_{X,\eta} \to \calo_T
\whotimes \A$, and since $\overline{\Phi}$ is symplectic, its
codifferential must be surjective. By Nakayama Lemma, this means
that both $\overline{\Phi}$ and $\Phi$ are automatically
isomorphisms. Therefore the correspondence $\Phi \mapsto \langle
\Phi(x_1),\dots,\Phi(x_n), \Phi(y_1), \dots, \Phi(y_n)\rangle$
identifies the set $\M_q(T)$ with the functor represented by a
closed subscheme in the total space of the $2n$-fold sum
$\left(J^\infty\D\right)^{\oplus 2n}$ of the bundle
$J^\infty\D$. This is our representing scheme $\M_q$.

It remains to note that $\M_q$ is naturally a $\WDA$-torsor (to
check that $\M_q$ is not only flat over $X$ but faithfully flat, one
uses Lemma~\ref{drboux}). Moreover, setting $\Phi \mapsto
\overline{\Phi}$ gives a natural map $\M_q \to \M_s$ compatible with
the map $\WDA \to \WSA$.
\endproof

\begin{remark} We note that the equivalence between torsors and
quantizations given in Lemma~\ref{qdrb} in fact goes through objects
of a third type: quantum-type deformations of the jet bundle
$J^\infty\calo_X$ in the tensor category of pro-vector bundles on
$X$ equiped with a flat connection. This might be useful, for
instance, in comparing our approach with that of A. Yekutieli -- any
isomorphism between two jet-bundle deformations by definition
induces a gauge equivalence in the sense of \cite{Y}, so that the
jet bundle deformations by definition satisfy the local differential
triviality condition of \cite{Y}.
\end{remark}

To make use of Lemma~\ref{qdrb}, we need some information of the
structure of the Lie algebra $\Der D$. Recall that every derivation
$d \in \Der D$ is almost inner -- namely, it can be obtained as the
commutator with an element
$$
\wt{d} \in h^{-1}D \subset D \otimes_{k[[h]]} k((h)).
$$ 
The vector space $h^{-1}D$ is closed under the commutator bracket
and forms a Lie algebra. Denote this Lie algebra by $\G$. Its center
coincides with the scalars $h^{-1}k[[h]] \cdot 1 \subset h^{-1}D$,
and we have a central extension of Lie algebras
\begin{equation}\label{cent}
\begin{CD}
0 @>>> k[[h]] @>{h^{-1}}>> \G @>>> \Der D @>>> 0.
\end{CD}
\end{equation}
For every $p \geq 0$, let $(\Der D)_{>p} = h^{p+1} \Der D \subset
\Der D$ be the subspace of $d \in \Der D$ such that for all $a \in
D$, $d(a)=0 \mod h^{p+1}$. Then $(\Der D)_{>p} \subset \Der D$ is a
Lie algebra ideal. Denote the quotient $\Der D/(\Der D)_{>p}$ by
$(\Der D)_p$. For every $p$, denote
$$
\G_p = \G/h^pD = \G/h^{p+1}\G.
$$
The extension \eqref{cent} is compatible with these quotients ---
for every $p \geq 0$, the Lie algebra $\G_p$ is a central extension
of the Lie algebra $(\Der D)_p$ by the vector space
$k[h]/h^{p+1}$. The kernel of the surjective map $\G_{p+1} \onto
\G_p$ is the space $\A = D/hD$ of functions on the standard
symplectic polydisc. The kernel of the map $(\Der D)_{p+1} \onto
(\Der D)_p$ is the vector space $\HH = \A/k \cdot 1$ of Hamiltonian
vector fields on the polydisc. All in all, for every $p \geq 0$ we
have a commutative diagram of the following type:
\begin{equation}\label{commu}
\begin{CD}
@.           0       @.       0      @.         0        @.   @.  \\
@.         @VVV             @VVV              @VVV            @.  \\
0  @>>> h^p \cdot k  @>>> h^p\cdot\A @>>>  h^p \cdot \HH   @>>> 0   \\
@.        @V{h}VV             @VVV              @VVV            @.  \\
0  @>>> k[h]/h^{p+2} @>>>  \G_{p+1}   @>>> (\Der D)_{p+1} @>>> 0   \\
@.         @VVV             @VVV              @VVV            @.  \\
0  @>>> k[h]/h^{p+1} @>>>    \G_p     @>>> (\Der D)_p     @>>> 0   \\
@.         @VVV             @VVV              @VVV            @.  \\
@.           0       @.       0      @.         0        @.   @.
\end{CD}
\end{equation}
The second and the third row, as well as the second and the third
column are extensions of Lie algebras. The first row and the first
column are extensions of the Lie algebra modules. Moreover, the
modules and the extension in the first column are trivial (and the
Lie algebra extensions in rows two and three are central).

Denote by $p$ the composition map $p:\Aut D \to \Sympl \A \to
Sp(2d)$; this is a projection from $\Aut D$ onto the symplectic
group. Its kernel $\Ker p \subset \Aut D$ is unipotent, and the
projection itself admits a splitting -- we have a semi-direct
product decomposition
$$
\Aut D \cong \Ker p \rtimes Sp(2d).
$$
This is in fact just the Levi splitting; to construct it explicitly,
one interprets the Weyl algebra $D$ as the completed universal
enveloping algebra of the Heisenberg Lie algebra. The group $Sp(2d)$
obviously acts on the Heisenberg algebra, hence on $D$. Moreover,
since $Sp(2d)$ is semisimple, the central extension \eqref{cent}
splits when restricted to $Sp(2d) \subset \Aut D$, so that the
semi-direct product decomposition lifts to the central
extension. This immediately allows us to integrate this central
extension to an extension of Harish-Chandra pairs. Indeed, it
suffices to integrate the extension of the kernel $\Ker p \subset
\Aut D$ and take the semi-direct product with $Sp(2d)$; since $\Ker
p$ is unipotent, constructing the extension of the group $\Ker p$
presents no problems (the center of this extension is the vector
space $k[[h]]$ considered as an additive group). By abuse of
notation, denote the whole extended Harish-Chandra pair by $\G$
(earlier this was used to denote its Lie algebra part). We will also
need to consider the quotient $\G/k$, where $k \subset k[[h]]$ is
the $1$-dimensional subspace of scalars; this quotient will be
denoted by $\Gb = \G/k$.

By construction we have a surjective map $\Gb \to \WSA$; denote its
kernel by $D^\dg$. The Lie algebra of the group $D^\dg$ is $D$
itself with the commutator bracket. The group $D^\dg$ is unipotent;
it is the product of the additive group $k$ and the group $D^u$ of
invertible elements in $D$ of the form $1+f$, where $f$ lies in the
ideal $1 + \m_\A+ hD \subset D$.

This integration procedure is compatible with the filtration by
degrees of $h$, so that we obtain the quotient Harish-Chandra pairs
$\WDAp$, $\G_p$, and $\Gb_p$, $p \geq 0$.

The extensions in the right-hand side column of \eqref{commu} are
non-trivial, both on the level of groups and on the level of Lie
algebras, with the following important exception.

\begin{lemma}\label{D1}
The Harish-Chandra extension
$$
\begin{CD}
0 @>>> \HH @>>> (\Der D)_1 @>>> (\Der D_0) @>>> 0
\end{CD}
$$
splits into a semidirect product, $(\Der D)_1 \cong \HH \rtimes
(\Der D)_0$.
\end{lemma}

\proof{} Notice that the Weyl algebra $D$ has a canonical
antiinvolution $\iota$ defined by $\iota(x_i)=x_i$,
$\iota(y_i)=y_i$, $\iota(h)=-h$. Consider the algebra $D_1 =
D/h^2D$. This is an associative algebra. Decomposing it with respect
to the eignevalues on $\iota$, we obtain a canonical
$\iota$-equivariant vector space identification $D_1 \cong \A \oplus
h \A$, $\iota(h)=-h$ and $\iota = \id$ on $\A$. Since $\iota$ is an
antiinvolution, under this identification the product in $D_1$ is
equal to
\begin{equation}\label{star}
a * b = ab + h\{a,b\},
\end{equation}
where $\{-,-\}$ is the Poisson bracket in $\A$. Moreover, the
Poisson bracket extends to a Lie bracket in $D_1$ defined
by 
$$
\{a,b\} = \frac{1}{h}\left(\wt{a}\wt{b}-\wt{b}\wt{a}\right) \mod h^2
$$ 
for any $a,b \in D_1$ lifted to $\wt{a},\wt{b} \in D$. Again, since
$\iota$ is an antiinvolution, the bracket in $D_1$ coincides with
the Poisson bracket on $\A$ extended to $\A \oplus h\A \cong
\A[h]/h^2$ in the $h$-linear way.

Denote by $(Der D)_1'$ the Lie algebra of derivations of the algebra
$D_1$ which are also derivations with respect to the Lie bracket. We
have a natural projection $(\Der D)'_1 \to (\Der D)_0$, and its
kernel is isomorphic to the vector space $\HH$. By definition, all
derivation of the algebra $D_1$ which come from derivations of the
algebra $D$ do preserve the Lie bracket, so that we have a natural
injective map $(\Der D)_1 \to (\Der D)'_1$. 

Since both $(\Der D)_1$ and $(\Der D)_1'$ project to $(\Der D)_0$
with the same kernel, an injective map $(\Der D)_1 \to (\Der D)'_1$
must be an isomorphism. Thus to prove the Lemma, it suffices to
construct a splitting $(\Der D)_0 \to (\Der D)_1'$. This is
immediate: the algebra $(\Der D)_0 \cong \HH$ acts naturally on $D_1
\cong \A \oplus h\A$, and this action preserves both the bracket and
the product \eqref{star}.
\endproof

The same argument also provides a splitting on the level of groups
and on the level of Harish-Chandra pairs. It is also easy to see
that this splitting extends to a splitting of the projection $\Gb_1
\to \Gb_0$.

\subsection{Categorical quantization.}\label{cate.sub}
At a suggestion of V. Drinfeld, we conclude this section with some
remarks on categorical aspects of quantization (this will not be
used in the rest of the paper; similar results were obtained earlier
in the analytic setting by P. Polesello and P. Schapira
\cite{PS}). Note that since the ideal $\m_A + hD \subset D$ is
invariant with respect to the group part $\Aut D$ of the
Harish-Chandra pair $\WDA$, the projection $\D^\dg \to k=D^\dg/D^u$
is also preserved by $\Aut D$. Therefore the group part $\Gb_{gr}$ of
the Harish-Chandra pair $\Gb$ admits a canonical product
decomposition
$$
\Gb_{gr} \cong k \times \left(\Gb_{gr}/k\right).
$$ 
Then starting from $\Gb$, we can define a different Harish-Chandra
pair $\Gb'$: it has the same Lie algebra part, but the group part is
replaced with the product $\Gb'_{gr} = k^* \times (\Gb_{gr}/k)$. In
other words, $\Gb'_{gr}$ is the extension of the group $\Aut D$ by
the group $D^*$ of invertible elements in the algebra $D$.

\begin{defn}\label{inte.defn}
A quantization $\D$ of an $S$-manifold $X$ is called {\em integral}
if the corresponding $\WDA$-torsor $\M_q$ lifts to a torsor $\M_q'$
over the Harish-Chandra pair $\Gb'$.
\end{defn}

\begin{remark}
We will see (in the end of Section~\ref{per}) that a canonical
quantization in the sense of Definition~\ref{can.defn} is always
integral; so, integral quantizations do exist, at least for
manifolds which admit a canonical quantization (such as, for
instance, manifolds with $H^i(X,\calo_X)=0$ for $i=1,2$).
\end{remark}

Assume given a symplectic $S$-manifold $X$ and an integral
quantization $\D$ of $X$, and denote by $\M_q'$ the lifting of the
corresponding torsor $\M_q$ to a $\Gb'$-torsor.  By construction, the
sheaf $\D$ is the sheaf of flat sections of the algebra bundle
$J^\infty\D$, so that the category of sheaves of left $\D$-modules
is equivalent to the category of pro-coherent sheaves of
$\calo_X$-modules equiped with a structure of a left module over
$J^\infty\D$ and a compatible flat connection. But it easy to see
that the torsor $\M'_q$ is transitive in the sense of
Definition~\ref{trans}. Therefore by Lemma~\ref{desce} the latter
category is in turn equivalent to the category of $\Gb'$-equivariant
pro-coherent sheaves of $\calo_{\M'_q}$-modules on $\M'_q$ equiped
with a $\Gb'$-equivariant structure of a left module over the Weyl
algebra $D$.

However, it turns out that there exists a different description of
this category which does not use the torsor $\M'_q$, nor
the torsor $\M_q$ (nor indeed the quantization $\D$). Namely, we
have the following.

\begin{prop}\label{cate}
Let $\D$ be a canonical quantization of a symplectic $S$-manifold
$X$. Let $\M_s$ be the torsor of formal symplectic coordinate
systems of $X$, and let the Harish-Chandra pair $\Gb'$ act on $\M_s$
through the quotient map $\Gb' \to \WSA$. Then the category of
sheaves of left $\D$-modules on $X$ is equivalent to the category of
$\Gb'$-equivariant pro-coherent sheaves $\E$ of
$\calo_{\M_s}$-modules on $\M_s$ such that the action of the kernel
$D^* \subset \Gb'$ of the projection $\Gb' \to \WSA$ on $\E$ extends
to a structure of a left module over $D$.
\end{prop}

To establish this equivalence, consider the natural action of the
Harish-Chandra pair $\Gb'$ on $D^*$, and let
$$
\Gt = D^* \rtimes \Gb'
$$ 
be the semi-direct product. We have a natural projection $\Gt \to
\WSA$, and its kernel is the semidirect product of the group $D^*$
with itself, with action by conjugation; this kernel is of course
canonically isomorphic to the product $D^* \times D^*$, so that
$\Gt$ is an extension of $\WSA$ by $D^* \times D^*$. We denote the
left and right copies of $D^* \subset D^* \times D^* \subset \Gt$ by
$D^*_l$ and $D^*_r$.

We have a natural embedding $\Gb' \to \Gt$ (it restricts to the
diagonal embedding $D^* \to D^* \times D^*$ on $D^* \subset \Gb'$);
we also have two projections $\tau_L,\tau_R:\Gt \to \Gb'$, whose
kernels are $D^*_r, D^*_l \subset D^* \times D^* \subset \Gt$. The
natural $\Gb'$-action on the formal Weyl algebra $D$ naturally
extends to a $\Gt$-module structure such that $D^* \times D^*
\subset \Gt$ acts by multiplication on the left and on the right. We
will denote $D$ with this $\Gt$-module structure by $D_o$. We also
need to consider two other $\Gt$-module structures on $D$: $D_L$
will be $D$ on which $\Gt$ acts through the projection $\tau_L:\Gt
\to \Gb'$, and $D_R$ will be $D$ on which $\Gt$ acts through
$\tau_R:\Gt \to \Gb'$. The $\Gb'$-action on $D$ is always the standard
one. Moreover, the action of $\Gt$ on $D_L$, $D_R$ and $D_o$ always
gives the standard action when restricted to $\Gb' \subset \Gt$.

Finally, we note that the $\Gt$-module structure on $D_R$ and $D_L$
is compatible with the algebra structure in $D$. The $\Gt$-module
$D_o$ is a $D_L-D_R$ bimodule, $D_L$ acts on the left, $D_R$ acts on
the right, both actions compatible with $\Gt$.

\proof[Proof of Proposition~\ref{cate}.] As noted above, the
category of sheaves of left $\D$-modules on $X$ is equivalent to the
category of $\Gb'$-equivariant pro-coherent sheaves of
$\calo_{\M'_q}$-modules on the torsor $\M'_q$ equiped with a
$\Gb'$-equivariant structure of a left module over the Weyl algebra
$D$. Let the Harish-Chandra pair $\Gt$ act on $\M'_q$ through the
projection $\tau_R:\Gt \to \Gb'$, so that the right-hand copy $D^*_r
\subset \D^* \times D^* \subset \Gt$ acts trivially. Then this
category is tautologically the same as the category $\CC_0$ of
$\Gt$-equivariant pro-coherent sheaves of $\calo_{\M'_q}$-modules on
$\M'_q$ on which $D^*_r$ acts trivially, and equiped with a
$\Gt$-equivariant structure of left $D_R$-module.

Taking tensor product over $D_R$ with the bimodule $D_o$ sends left
$D_R$-modules into left $D_L$-modules, and it is immediate to check
that left $D_R$-modules on which $D^*_l$ acts trivially are sent
into left $D_L$-modules on which $D^*_l$ acts by left
multiplication. This tensor product functor is obviously an
equivalence. We conclude that the category $\CC_0$ is naturally
equivalent to the category $\CC_1$ of $\Gt$-equivariant pro-coherent
sheaves of $\calo_{\M'_q}$-modules on $\M'_q$ equiped with a
$\Gt$-equivariant structure of a left $D_L$-module, such that $D_l^*
\subset \Gt$ acts by left multiplication.

Consider now the right-hand copy $D^*_r \subset D^* \times D^*
\subset \Gt$. This subgroup acts freely on $\M'_q$, and we have
$\M'_q/D^* \cong \M_s$, so that $\M'_q$ is a principal $D^*$-bundle
over $\M_s$. On the other hand, $D_r^*$ acts trivially on
$D_L$. Applying descent with respect to $D^*_r$, we identify $\CC_1$
with the category of $\Gb' = \Gt/D^*_r$-equivariant pro-coherent
sheaves of $\calo_{\M_s}$-modules on $\M_s$ equiped with a
structure of a left $D$-module such that $D^* \subset \Gb'$ acts by
left multiplication. Since the $D$-module structure is uniquely
defined by the $D^*$-action, this proves the Proposition.

We note that one can also consider the bundles $\D_R$, $\D_L$ and
$\D_o$ on $\M_s$ obtained by descent from the $\Gt$-modules $D_R$,
$D_L$ and $D_o$.  They are naturally equivariant with respect to $\Gb'
= \Gt/D^*_r$; $\D_R$ and $\D_L$ are bundles of algebras, while
$\D_o$ is a $\D_L-\D_R$-bimodule. The bundle $\D_L$ is constant, and
the group $D^* \subset \Gb'$ acts on it by conjugation. The bundle
$\D_R$ is not constant, but the group $D^* \subset \Gb'$ acts on it
trivially. The categories of $\Gb'$-equivariant left $\D_R$ and
$\D_L$-modules are equivalent, and the equivalence is given by
tensoring with $\D_o$.
\endproof

The reason Proposition~\ref{cate} is interesting is that the
right-hand side of the established equivalence is defined {\em a
priori}, without any reference either to $\D$ or to $\M_q$. Thus it
gives a perfectly well-defined abelian category $Quan(X)$ for {\em
any} symplectic $S$-manifold $X$, not only for an admissible one.
Loosely speaking, in terms of this category, the problem of finding
an integral quantization of a symplectic $S$-manifold $X$ becomes
the problem of finding an object of rank 1 in $Quan(X)$ (the object
$\D_o$ described above). Given such an object, the
$\WSA$-equivariant algebra $\D_R$ on $\M_s$ is recovered as the
endomorphism algebra of $\D_o$, and the quantization $\D$ is
obtained by localizing $\D_R$ with respect to $\WSA$. This gives an
alternative approach to the quantization problem (which is however
pretty close in essense, if not in form, to our Lemma~\ref{qdrb}).

\begin{remark} As pretty much everything in this paper, the argument
in this subsection is not terribly original. In this particular
case, we learned the idea from V. Drinfeld; its origins are
attributed to J. Bernstein, P. Deligne and M. Kontsevich. The
influence of Kontsevich of course looms large over the whole
subject, although in this text we have deliberately used
old-fashioned arguments independent from the Formality Theorem.
Notice also that \cite{PS} is based on earlier work on contact
manifolds done by the Japanese school, see e.g. \cite[\S 8.2]{K} or
\cite[Chapter 7]{KS}.
\end{remark}

\begin{remark} If one uses the existence of canonical quantization
for affine manifolds (which follows from Theorem~\ref{main}, since
for an affine $X$ we have $H^i(X,\calo_X)=0$), then one can recast
the construction of the category $Quan(X)$ in a \v{C}ech-cocycle
style. Namely, the quantizations exist locally; being canonical,
they are isomorphic on intersections. The compatiblity isomorphisms
do not agree on triple intersections, but the difference between
them is an {\em inner} automorphism $\gamma_{i,j,k}$ of the algebra
$\D$. Moreover, if instead of just quantizations one considers pairs
of a quantization $\D$ and a lifting of the corresponding torsor
$\M_q$ to a $\Gb$-torsor, then the automorphisms $\gamma_{i,j,k}$
come equiped with a lifting to an element of $\D^*$. One then uses
these elements to glue together the local categories by a standard
construction. The argument in this subsection is essentially the
same, but a \v{C}ech covering is replaced with the torsor $\M_s$;
this allows to avoid using Theorem~\ref{main}.
\end{remark}

\section{The non-commutative period map.}\label{per}

We can now define the non-commutative period map and prove
Theorem~\ref{main}. Fix an $S$-manifold $X$ of dimension $n=2d$
equiped with a symplectic form $\Omega$. Let $\M_s$ be the bundle
of symplectic formal coordinate systems on $X$. By Lemma~\ref{qdrb},
the set $Q(X,\Omega)$ of quantizations of $\langle X, \Omega
\rangle$ is in one-to-one correspondence with the set
$H^1_{\M_s}(X,\WDA)$ of the liftings of $\M_s$ to a $\WDA$-torsor
$\M_q$.

Recall that the central extension \eqref{cent} is integrated to a
central extension
$$
\begin{CD}
1 @>>> k[[h]] @>>> \G @>>> \WDA @>>> 1
\end{CD}
$$
of Harish-Chandra pairs. For any $\WDA$-torsor $\M_q$, the
localization of the trivial $\WDA$-module $k[[h]]$ is the constant
local system $\calo_X[[h]]$ on $X$. By Proposition~\ref{nab}, we
have an obstruction map from the set $Q(X,\Omega) \cong
H^1_{\M_s}(X,\WDA)$ to the second cohomology group
$H^2_{DR}(X)[[h]]$ -- it sends a torsor $\M$ to the class which
obstructs the lifting of $\M$ to a torsor over $\G$.

\begin{defn}\label{per.defn}
The obstruction map
$$
\Per:Q(X,\Omega) \to H^2_{DR}(X,A[[h]]) \cong H^2_{DR}(X)[[h]]
$$
associated to the extension \eqref{cent} is called the {\em
non-commutative period map}.
\end{defn}

\proof[Proof of Theorem~\ref{main}.]  For every $p \geq 0$,
integrate the central extension \eqref{cent} to a Harish-Chandra
extension of the quotient Harish-Chandra pair $\WDAp$ by the trivial
module $h^{-1}k[h]/h^pk[h]$. By abuse of notation, denote the whole
extended Harish-Chandra pair by $\G_p$. The commutative diagram
\eqref{commu} can be considered as a diagram of Harish-Chandra pairs
and their extensions.

We have to prove that the period map $\Per:Q(X,\Omega) \to
H^2_{DR}(X)[[h]]$ is injective, maps any quantization to a power
series with constant term $[\Omega]$, and that any splitting
$F:H^2_{DR}(X) \to H^2_F(X)$ induces an isomorphism
$$
Q(X,\Omega) \cong F([\Omega]) + hH^2_F(X)[[h]].
$$
Fix such a splitting $F:H^2_{DR}(X) \to H^2_F(X)$. To simplify
notation, denote the set $H^1_{\M_s}(X,\WDAp)$ by $Q_p$, and denote
by $\Per_p$ the obstruction map $Q_p \to
H^2_{DR}(X,A[h]/h^p)$. Since $(\Der D)_0$ is simply the algebra
$\HH$ of Hamiltonian vector fields, the set $Q_0$ consists of one
point, namely, the $\WSA$-torsor $\M_s$. By the remarks after
Lemma~\ref{drb}, we have
$$
\Per_0(Q_0) = \loc(\M_s,[\omega]) = [\Omega] \in H^2_{DR}(X).
$$
By induction, it suffices to prove that for every $l > 0$, the map
$\Per_l$ is injective, and the projection $F$ identifies its image
with $H^2_F(X) \otimes_k k[h]/h^{l+1}$. We may assume the claim
proved for all $l \leq p$ and consider the case $l = p+1$. Moreover,
we may fix a torsor $\M \in Q_p$. Once we do it, it suffices to
prove that the period map $\Per_{p+1}$ is injective on the set
$$
H^1_\M(X,\WDApt),
$$ 
and that it sends this set to a torsor over $H^2_F(X) \cdot h^{p+1}
\subset H^2_{DR}(X)[h]/h^{p+2})$.

By \eqref{commu}, the Harish-Chandra pair $\G_{p+1}$ is an extension
of the Harish-Chandra pair $(\Der D)_p$ by the module
$$
V = \left(k[h]/h^{p+2} \oplus \A \cdot h^p\right)/k \cdot h^p.
$$
Consider the submodule $U = k[h]/h^{p+2} \subset V$, and denote the
quotient module by $W = V/U = \A/k = \HH$. Thus we have a short
exact sequence \eqref{ab} and a Harish-Chandra extension of the type
considered in Lemma~\ref{oct}, with $\GH = \WDAp$. The intermediate
Harish-Chandra extension $\GHt$ is given by $\WDApt$. As in
Lemma~\ref{oct}, denote the localizations of the $\WDAp$-modules
$U$, $V$ and $W$ with respect to the torsor $\M$ by $\U$, $\V$ and
$\W$.

\begin{lemma}\label{unobs}\mbox{}
Assume that the manifold $X$ is admissible.
\begin{enumerate}
\item The canonical map $H^2_{DR}(X,\V) \to H^2_{DR}(X,\W)$ is
trivial. 
\item The canonical map $H^1_{DR}(X,\W) \to H^2_{DR}(X,\U)$ is
injective.
\end{enumerate}
\end{lemma}

\proof{} By the long exact sequence associated to \eqref{ab},
\thetag{ii} is equivalent to saying that the canonical map
$H^1_{DR}(X,\V) \to H^1_{DR}(X,\W)$ is trivial. In other words, we
have to prove that the map $H^l_{DR}(X,\V) \to H^l_{DR}(X,\W)$ is
trivial for $l=1,2$.

Note that the $\WDAp$-module structure on $U$, $V$ and $W$ is
obtained by restriction from a $\WSA$-module structure by the
canonical map $\WDAp \to \WSA$. Therefore the localizations $\U$,
$\V$ and $\W$ do not depend on the torsor $\M$. Moreover, since
$k[h]/h^{p+2} \cong k \cdot h^{p+1} \oplus k[h]/h^{p+1}$ as vector
spaces, we have $V \cong \U \oplus k[h]/h^{p+1}$, and the map $V \to
W$ is trivial on the second summand. Therefore it suffices to prove
that the surjection $\A \to \HH \cong W$ induces a trivial map
$$
H^l_{DR}(X,\loc(\M_s,\A)) \to H^l_{DR}(X,\W) \cong
H^l_{DR}(X,\loc(\M_s,\HH))
$$
for $l=1,2$. Since $\HH = \A/k$, this is in turn equivalent to
saying that the map
$$
H^l_{DR}(X) \cong H^l_{DR}(X,\loc(\M_s,k)) \to
H^l_{DR}(X,\loc(\M_s,\A))
$$
is surjective for $l=1,2$. But we know that $\loc(\M_s,\A) \cong
J^\infty\calo_X$. Therefore $H^l_{DR}(X,\loc(\M_s,\A)) \cong
H^l(X,\calo_X)$, and the claim becomes the definition of
admissibility.
\endproof

Let $c \in H^2_{DR}(X,\V)$ be the obstruction class associated to
the torsor $\M$ and the extension $\G_{p+1}$. By Lemma~\ref{unobs},
the class $c$ restricts to zero in $H^2_{DR}(X,\W)$. Thus the
assumptions of Lemma~\ref{oct} are satisfied. We conclude that the
period map
\begin{multline*}
\Per_{p+1}:H^1_{\M}(X,\WDApt) \to\\
\to H^2_{DR}(X,\U) = H^2_{DR}(X,A[h]/h^{p+2})
\end{multline*}
is compatible with the $H^1_{DR}(X,\W)$-action on both sides. By
Proposition~\ref{nab}, the left-hand side is a
$H^1_{DR}(X,\W)$-torsor. Moreover, we have $H^1_{DR}(X,\W) \cong
H^2_F(X)$. Thus to prove the inductive step and the Theorem, it
suffices to prove that the $H^1_{DR}(X,\W)$-action on the right-hand
side is free. This is exactly Lemma~\ref{unobs}~\thetag{ii}.
\endproof

We finish this Section with some observations on the canonical and
integral quantizations. By our definition of the non-abelian period
map, a symplectic $S$-manifold $X$ is canonical in the sense of
Definition~\ref{can.defn} if and only if the corresponding
$\WDA$-torsor $\M_q$ over $X$ lifts to a $\Gb$-torsor. One can try
to construct such a lifting by going step-by-step though the
quotient groups $\Gb_p$. Due to Lemma~\ref{D1}, there are no
obstructions to this at the first step -- the symplectic coordinate
torsor $\M_s$ always lifts to a $\Gb_1$-torsor, and a preferred
lifting is given by the torsor $\M'_1$ induced by means of the
splitting map $\WSA \cong \Gb_0 \to \Gb_1$. There may be other
liftings, but by Theorem~\ref{main} they all become isomorphic after
taking quotient by the center $h \cdot k \subset \Gb_1$.

At every further step, there is an obstruction class lying in the
group $H^2(X,\calo_X)$. However, we note that these obstructions do
not depend on the entire $\Gb_1$-torsor $\M'_1$, but only on its
reduction $\M_1 = \M'_1/k$ to $\Gb_1/k$. Indeed, we have
$$
\Gb \cong \Gb_1 \times_{\Gb_1/k} \Gb/(h \cdot k),
$$
where $h \cdot k \subset hk[[h]]$ lies in the center $hk[[h]]$ of
the Harish-Cahndra pair $\Gb$. Then every $\Gb$-torsor $\M'_q$
is of the form
$$
\M'_q \cong \M'_1 \times_{\M_1} \Mb_q,
$$
where $\M'_1=\M_q/\Gb_{>1}$ is the reduction of $\M_q'$ to $\Gb_1$,
$\Mb_q= \M'_q/(h \cdot k)$ is its reduction to $\Gb/(h \cdot k)$,
and $\M_1=\M'_1/k=\Mb_q/\Gb_{>1}$ is their common reduction to
$\Gb/((h \cdot k) \times \Gb_{>1}$. This means that in order to lift
a $\WDA$-torsor $\M_q$ to $\Gb$, it is sufficient, and necessary, to
lift it to $\Gb/(h \cdot k)$, and to lift its reduction $\M_1 =
\M_q/\Gb_{>1}$ to $\Gb_1$.

Analogously, one can try to construct step-by-step an integral
quantization in the sense of Definition~\ref{inte.defn}. Again,
there is no obstruction at the first step, and there is a preferred
lifting given by induction with respect to the splitting map. We
leave it to the reader to check that all other liftings are
parametrized by line bundles $L$ on $X$, and the leading two terms
of the period map for an integral quantization is of the form
$[\Omega] + hc_1(L)$, where $c_1(L)$ is the first Chern class. To
construct an integral quantization, it suffices to have a
$\WDA$-tors $\M_q$ which lifts to $\Gb'/k^*$, and whose reduction
$\M_1/\Gb'_{>1}$ to $\Gb'_1/k^*$ lifts to $\Gb'_1$.

In particular, if we are given a canonical quantization, then the
corresponding $\WDA$-torsor $\M_q$ comes equiped with a lifting to
$\Gb'/k^*$ -- indeed, by definition we have $\Gb'/k^* \cong \Gb/(h
\cdot k)$ and $\Gb'_1/k^* \cong \Gb_1/k$. Moreover, the reduction
$\M_1=\M_q/\Gb_{>1}$ is the preferred lifting of the symplectic
coordinate torsor $\M_s$, so that it lifts to $\Gb'_1$. We conclude
that a canonical quantization is integral in the sense of
Definition~\ref{inte.defn}.

\begin{remark}\label{rw}
One can actually compute the first of the obstruction classes,
namely, the obstruction to lifting a $\Gb_1$-torsor to $\Gb_2$. This
amounts to computing the cocycle that represents the extension
$\Gb_2 \to \G_1$; we will not do it here to save space, and only
state the answer -- up to a non-zero constant, the obstruction is
given by the so-called Rozansky-Witten class corresponding to the
trivalent graph with two vertices connected by three edges (for
Rozansky-Witten classes, see e.g. \cite{rowi}; for the proof in the
holomorphic setting, see in \cite{NT}).
\end{remark}

\section{Non-abelian cohomology.}\label{na}

We now return to the basics and prove the results announced in
Subsection~\ref{nas}. In the case of ordinary torsors, even in a
very general topos, everything is completely standard (\cite{Gi},
or, for example, a much shorter and nicer exposition in
\cite{Ga}). Unfortunately, we need to work with flat
connections. One can probably obtain all the results for free by
passing to the crystalline topos, but this raises the amount of high
science used to a completely disproportionate degree. For the
convenience of the reader, and for our own peace of mind, we will
give a proof of all the facts we need in down-to-earth terms. To
save space, the more standard parts of the proofs are left to the
reader.

\subsection{Linear algebra.}\label{la}
Recall that for any two objects $A$, $B$ in a fixed abelian
category, we can form the {\em extension groupoid} $\ext^1(B,A)$
whose objects are short exact sequences
$$
\begin{CD}
0 @>>> A @>>> \bullet @>>> B @>>> 0,
\end{CD}
$$
and whose morphisms are isomorphisms of the exact sequences
identical on $A$ and on $B$. The set of isomorphism classes of
objects in the groupoid $\ext^1(B,A)$ is the first $\Ext$-group
$\Ext^1(B,A)$. The groupoid $\ext^1(B,A)$ has an additional
structure of a symmetric monoidal category: the sum is given by the
Baer sum of extensions.

Fix objects $A$, $B$, and let $c \in \Ext^2(B,A)$ be an element in
the second $\Ext$-group. Represent $c$ by a four-term exact
sequence, Yoneda-style
\begin{equation}\label{4t}
\begin{CD}
0 @>>> A @>>> E_1 @>>> E_2 @>>> B @>>> 0.
\end{CD}
\end{equation}
In other words, we have a two-term complex $E_1 \to E_2$ whose
cohomology objects are $A$ and $B$. Recall that $c=0$ if and only if
there exist a complex
\begin{equation}\label{mon}
\begin{CD}
A @>{a}>> E @>{b}>> B
\end{CD}
\end{equation}
such that $\Ker b \cong E_1$, $\Coker a \cong E_2$, and the natural
sequence
$$
\begin{CD}
0 @>>> A @>>> \Ker b @>>> \Coker a @>>> B @>>> 0
\end{CD}
$$
is exact and isomorphic to \eqref{4t}. Diagrams of the form
\eqref{mon} form a groupoid, which we will denote by $\spl(c)$ (we
require maps between diagrams to be identical on $\Ker b \cong E_1$
and $\Coker a \cong E_2$). The groupoid $\spl(c)$ is naturally a
gerb over the symmetric monoidal groupoid $\ext^1(B,A)$ -- this
means that we have a sum functor $\ext^1(B,A) \times \spl(c) \to
\spl(c)$, a difference functor $\spl(c) \times \spl(c) \to
\ext^1(B,A)$, and natural compatibility morphisms between these
functors which turn $\spl(c)$ into a ``torsor'' over $\ext^1(B,A)$
in the obvious sense. Both the sum and the difference functor are
again given by the Baer sum construction.

This construction is functorial in the following way: every exact
functor $F$ between abelian categories induces a functor $F:\spl(c)
\to \spl(F(c))$ between groupoids $\spl(c)$ and $\spl(F(c))$.

Taking a different Yoneda representation \eqref{4t} for the same
element $c \in \Ext^2(B,A)$ gives an equivalent groupoid
$\spl(c)$. To make this quite canonical, one has to consider all
possible representations and treat objects of $\spl(c)$ as certain
diagrams of sheaves on the category of these representations. This
is very beautiful but too technical to describe here, see
\cite{Gr}. For our purposes, it suffices to carry a fixed Yoneda
representation in all the constructions.

Of course, if the class $c$ is not trivial, the groupoid $\spl(c)$
is empty. But it is important to define it anyway.

Finally, in proving Lemma~\ref{oct} we will need to consider the
following situation. Assume that the object $A$ in \eqref{4t} is the
middle term of a short exact sequence
\begin{equation}\label{shrt}
\begin{CD}
0 @>>> A_0 @>{a}>> A @>{b}>> A_1 @>>> 0.
\end{CD}
\end{equation}
We have the canonical Yoneda representation $E_1/A_0 \to E_2$ of the
class $b(c) \in \Ext^2(B,A_1)$. Assume in addition that we have
$b(c) = 0$. Then we claim that every object $s \in \spl(b(c))$
canonically defines a class $c_0 \in \Ext^2(B,A_0)$ such that $c =
a(c_0) \in \Ext^2(B,A)$. Indeed, we can take the class represented
by the exact sequence
$$
\begin{CD}
0 @>>> A_0 @>>> E_1 @>>> E @>>> B @>>> 0,
\end{CD}
$$
where $A_1 \to E \to B$ represents the object $s \in \spl(b(c))$. If
we twist the object $s$ by an extension $e \in \Ext^1(B,A_1)$, then
the corresponding class $c_0$ is replaced by $c_0 + d(e)$, where
$d:\Ext^1(B,A_1) \to \Ext^2(B,A_0)$ is the differential in the long
exact sequence associated to \eqref{shrt}. To prove it, it suffices
to notice that the twisting is by definition done via Baer sum with
the sequence
$$
\begin{CD}
0 @>>> A_0 @>>> A @>>> F @>>> B @>>> 0,
\end{CD}
$$
the Yoneda product of the sequence \eqref{shrt} and the sequence $0
\to A_1 \to F \to B \to 0$ which represents $e \in \Ext^1(B,A_1)$.

\subsection{Extensions.}
Assume given a Harish-Chandra pair $\GH$ and a $\GH$-module $V$.
Consider $V$ as an (additive) algebraic group. Let $\vv$ be the
tautogolical Harish-Chandra pair $\vv$. The Harish-Chandra
cohomology groups $H^\hdot(\vv,V)$ coincide with the cohomology
$H^\hdot(V,V)$ of the group (or Lie algebra) $V$. In particular,
classes in the group $H^1(\vv,V)$ correspond to Lie algebra
derivations $d:V \to V$.  Denote by $\tau_V \in H^1(\vv,V)$ the
tautological class -- namely, the one corresponding to the identity
map $\id:V \to V$.

Fix a Harish-Chandra extension of the type \eqref{extt}, and
consider the Hochschild-Serre spectral sequence 
$$
H^\hdot\left(\GH,H^\hdot(V,V)\right) \Rightarrow
H^\hdot\left(\GHh,V\right) 
$$
which computes the cohomology groups $H^\hdot(\GHh,V)$. The
$E^2$-term of this sequence contains in particular the group
$H^1(\vv,V)$, and the differential gives a map
$$
d:H^1(\vv,V) \to H^2(\GH,V).
$$
Applying $d$ to the tautological class $\tau_V$ gives an element 
$$
c = d\tau_V \in H^2(\GH,V)
$$ 
canonically associated to the extension.

This class $c$ is of course just the usual $2$-cocycle known both in
the theory of algebraic groups and in the theory of Lie
algebras. Out of the myriad equivalent ways to construct it, this
particular one has the advantage of only using the Hochschild-Serre
spectral sequence. Therefore it generalizes to Harish-Chandra pairs
without any additional work. We record explicitly one degenerate
case: when the group $G$ is trivial, the cocycle
$$
c \in H^2(\GH,V) = H^2(\h,V) = \Hom(\Lambda^2\h,V)
$$
is just the commutator map in the Lie algebra $\h_1$.

For every map $f$ between Harish-Chandra pairs, denote by $f^*$ the
restriction functor on the categories of modules.  Fix a particular
Yoneda representation of the class $c \in H^2(\GH,V) = \Ext^2(k,V)$,
and consider the groupoid $\spl(c)$. Fix a splitting of the exact
sequence \eqref{4t} considered as a sequence of vector spaces.
Since the composition $\pi \circ \rho$ factors through the map
$\eta:H \to 1$, the fixed vector-space splitting defines a canonical
object in the groupoid $\spl(\rho^*\pi^*(c))$. This gives a
trivialization of the corresponding gerbe, that is, an equivalence
$$
\spl(\rho^*\pi^*(c)) \cong \ext^1(k, \rho^*\pi^*(V)).
$$
Note now that by construction, the class $\pi^*c$ is trivial.
Therefore the groupoid $\spl(\pi^*(c))$ is non-empty. Recall that it
is also a gerbe over the extension groupoid
$\ext^1(k,\pi^*(V))$. Say that an object $s \in \spl(\pi^*(c))$ is a
{\em good splitting} if $\rho^*s$ is the tautological extension
represented by the cohomology class $\tau_V$. Analyzing the
Hochschild-Serre spectral sequence, we see that, since $c =
d(\tau_V)$, good objects exist, and that pairs $\langle s, f:\rho^*
\cong \tau_V\rangle$ form a gerbe over $\ext^1(k,V)$. Since our goal
is not to construct a general theory of group extensions but rather,
to have a skeleton theory sufficient for applications to torsors, we
will simply ignore this ambiguity and fix a good splitting $s \in
\spl(\pi^*(c))$ for every Harish-Chandra extension \eqref{extt}.

\subsection{Torsors.}
Fix an $S$-manifold $X$. Assume given a Harish-Chandra extension
\eqref{extt} and a $\GH$-torsor $\M$ over $X$. Denote by $\V =
\loc(\M,V)$ the localization of the $\GH$-module $V$ with respect to
$\M$. Consider the groupoid $\hh^1_{\M}(X,\GHh)$ of $\GHh$-torsors
$\M_0$ over $X$ equiped with an isomorphism $\pi_*\M_0 \cong \M$.

Let $c(\M) \in H^2_{DR}(X,\V)$ be the localization of the cohomology
class $c$ with respect to the torsor $\M$. This class comes equiped
with a Yoneda representation (obtained by the localization of the
fixed Yoneda representation of the class $c$). Moreover, for any
torsor $\M_0 \in \hh^1_{\M}(X,\GHh)$, the localization
$\loc(\M_0,\pi^*c)$ canonically coincides with $c(\M)$. We can set
$$
\M_0 \mapsto \loc(\M_0,s)
$$
and obtain a functor $\Lin:\hh^1_{\M}(X,\GHh) \to \spl(c(\M))$. The
crucial part of both Proposition~\ref{nab} and Lemma~\ref{oct} is
the following fact.

\begin{lemma}\label{eqq}
Assume that the $\GH$-torsor $\M$ is transitive in the sense of
Definition~\ref{trans}. Then the functor
\begin{equation}\label{lin}
\Lin:\hh^1_{\M}(X,\GHh) \to \spl(c(\M))
\end{equation}
is an equivalence of categories.
\end{lemma}

\proof{} We will define an inverse equivalence. Assume given a
splitting $s_X \in \spl(c(\M))$. Consider the fixed Yoneda
representation
$$
\begin{CD}
0 @>>> \V @>>> \E_1 @>>> \E_2 @>>> \calo_X @>>> 0
\end{CD}
$$
of the class $c(\M)$ split by $s_X$. Denote by $\sigma:\M \to X$ the
projection. By definition of the localization functor, the diagram
$$
\begin{CD}
0 @>>> \sigma^*\V @>>> \sigma^*\E_1 @>>> \sigma^*\E_2 @>>>
\calo_{\M} @>>> 0
\end{CD}
$$ 
is a diagram of {\em constant} vector bundles on $\M$: it is
isomorphic to
$$
\begin{CD}
0 @>>> \rho^*V @>>> \rho^*E_1 @>>> \rho^*E_2 @>>> \rho^*k @>>> 0,
\end{CD}
$$
where $\rho:\M \to \Spec k$ is the projection to the point, and
\begin{equation}\label{vs}
\begin{CD}
0 @>>> V @>>> E_1 @>>> E_2 @>>> k @>>> 0,
\end{CD}
\end{equation}
is the diagram representing the class $c$ of the Harish-Chandra
extension $\GHh$. We have two splittings of this diagram: one is
given by the fixed $\GHh$-equivariant splitting $s$ of the diagram
\eqref{vs}, the other is given by $\sigma^*(s_X)$.

For every point $m \in \M$, denote by $A_m$ the set of all
isomorphisms
\begin{equation}\label{iso}
\phi:s_m \cong \sigma^*(s_X)_m
\end{equation}
between the fibers of these splittings at the point $m$ (which we
consider as splittings in the category of vector spaces). Denote by
$\M_s$ the set of pairs
$$
\langle m, \phi \in A_m \rangle
$$
of a point $m \in \M$ and an isomorphism $\phi \in A_m$. Each of the
sets $A_m$ is naturally a torsor over the vector space $\Hom(k,V) =
V$. Therefore $\M_s/\M$ is a torsor over the constant bundle
$\sigma^*\V$ on $\M$. This makes it into a scheme over $\M$, in fact
into a (pro)$S$-manifold.

Since the splitting $s$ is in fact a splitting in the category of
$\GHh$-modules, the group $G_1$ acts naturally on the $S$-manifold
$\M_s$: it acts on $m$ through the quotient $G = G_1/V$, and it acts
on $\phi$ by acting on the left-hand side of \eqref{iso}. Since the
stabilizer of a point $m \in \M$ is the subgroup $V \subset G_1$,
and the set $A_m$ is a $V$-torsor, the whole $\M_s$ is a
$G_1$-torsor. Moreover, the natural $G_1$-action naturally extends
to a transitive $\h_1$ action, and turns $\M_s$ into a well-defined
transitive $\GHh$-torsor $\M_s \in \hh^1_{\M}(X,\GHh)$.

\medskip

The correspondence $s_X \mapsto \M_s$ is obviously functorial in
$s_X$ and gives a functor $\spl(c(\M)) \to \hh^1_{\M}(X,\GHh)$.
This is the desired inverse equivalence, more or less by definition;
the proof is easy, and we leave it to the reader.
\endproof

This Lemma immediately yields Proposition~\ref{nab} and
Lemma~\ref{oct}. Indeed, it allows to rewrite both statements as
claims about the groupoid $\spl(c(\M))$, and these claims
immediately follow from the homological considerations of
Subsection~\ref{la}.

\section{Generalizations.}\label{ex}

Theorem~\ref{main} admits two immediate generalizations -- indeed,
both could have been incorporated directly into its statement, and
we did not do so only out of desire to keep the statement down to a
reasonable size. We record both here, together with a result on
comparison with symplectic deformation theory of \cite{KV}.

\subsection{Equivariant situation.} Let $X$ be an $S$-manifold
equiped with an action of a reductive group $G$. Note that $G$ acts
naturally on the de Rham cohomology $H^\hdot_{DR}(X)$ and on the
coherent cohomology $H^\hdot(X,\calo_X)$.

\begin{defn} The $S$-manifold $X$ equiped with the $G$-action is
called {\em admissible in the $G$-equivariant sense} if the
canonical map
$$
\left(H^i_{DR}(X)\right)^G \to \left(H^i(X,\calo_X)\right)^G
$$
between the $G$-invariant parts of the respective cohomology groups
is surjective for $i=1,2$.
\end{defn}

\begin{prop}\label{p1}
Let $X$ be a symplectic $S$-manifold equiped with a $G$-action
which presevres the symplectic form, and assume that $X$ is
admissible in the $G$-equivariant sense. Then $X$ has a
$G$-equivariant quantization.
\end{prop}

\proof{} The proof of Theorem~\ref{main} works without any changes,
save for adding ``$G$-equivariant'' in appropriate places. Note that
$G$-equivariant local systems should be understood in a ``stupid''
way -- as $G$-equivariant vector bundles equiped with a
$G$-invariant flat connection. In particular, genuine
$G$-equivariant cohomology groups $H_G^\hdot(X)$ do not enter into
the picture.
\endproof

The canonical quantization, being canonical, is equivariant with
respect to any possible group action. This allows to define and
construct quantizations of admissible global quotients by a finite
group -- indeed, for a finite group $G$, a quotient $X = Y/G$ is
admissible if and only if $Y$ is admissible in the $G$-equivariant
sense. Quantization of arbitrary admissible Deligne-Mumford stacks
is more delicate, and we prefer to postpone this investigation to a
future paper.

Another situation when Proposition~\ref{p1} might be useful is when
we want to quantize a symplectic manifold which is not admissible in
the sense of Definition~\ref{adm}. For example, given an
$S$-manifold $X$ with $H^1(X,\calo_X) = H^2(\calo_X)=0$ and a line
bundle $L$ on $X$, one can consider the total space $Z$ of the
associated $\gm$-torsor over $X$. Typically tensor powers $L^k$, $k
\in \Z$ of the bundle $L$ will have non-trivial cohomology groups,
so that $H^i(Z,\calo_Z)$ would be large and $Z$ would not have a
chance of being admissible. However, since
$$
\left(H^i(Z,\calo_Z)\right)^\gm = H^i(X,\calo_X) = 0, \qquad i=1,2,
$$
the manifold $Z$ is always admissible in the $\gm$-equivariant sense.

\subsection{Comparison with symplectic deformations.}

Theorem~\ref{main} holds literally, with the same proof, when either
$S$, or $X$, or both are allowed to be formal schemes -- indeed, all
we ever used of a scheme was a formal neighborhood of its closed
point. This allows for comparison with \cite{KV}. The main result of
\cite{KV} is the following.

\begin{theorem}[{{\cite[Theorem 1.1]{KV}}}]\label{kv}
Let $X$ be an admissible manifold over the field $k$. Assume that
$X$ is equiped with a nondegenerate symplectic form
$\Omega_0$. Then the pair $\langle X, \Omega \rangle$ admits a
universal formal deformation $\X/S$. Moreover, the cohomology class
$$
[\Omega] \in H^2_{DR}(\X/S) \cong H^2_{DR}(X) \otimes \calo_S
$$
of the relative symplectic form $\Omega \in \Omega^2(\X/S)$ defines
an embedding $S \to H^2_{DR}(X)$, and every splitting $H^2_{DR}(X)
\to H^2_F(X)$ of the natural embedding $H^2_F(X) \to H^2_{DR}(X)$
identifies $S$ with the formal completion of the affine space
$H^2_F(X)$ at the point $[\Omega_0] \in H^2_F(X)$.\endproof
\end{theorem}

In general, the universal deformation $\X/S$ exists only as a formal
scheme. The precise meaning of universality will not be important
for us, see \cite{KV}. What is important is that $\X$ is smooth and
symplectic over $S$, so that we can apply Theorem~\ref{main} and
construct a quantization of $\X/S$. Having done this, we obtain a
non-commutative multiparameter deformation $\D_S$ of the structure
sheaf of the symplectic manifold $X/\Spec k$. The base of this
deformation is
$$
\overline{S} = \Delta \times S \subset \Delta \times H^2_{DR}(X),
$$
where $\Delta = \Spf k[[h]]$ is the formal disc, and $H^2_{DR}(X)$
is considered as an affine space. For any section $s:\Delta \to
\overline{S}$ of the natural projection $\overline{S} \to \Delta$,
the pullback $s^*\D_S$ is a quantization of the manifold
$X$. Algebraically, every such section $s$ is given by a formal
power series $P_s \in H^2_{DR}(X)[[h]]$.

\begin{lemma}\label{basech}
Assume that the chosen quantization $\D_S$ of the $S$-manifold $\X$
is canonical. Then the non-commutative period map sends the
quantization $s^*\D_S$ to the formal power series $P_s$.
\end{lemma}

\proof{} This is immediate from the definitions. Indeed, since
$\D_S$ is the canonical quantization of $\X/S$, its non-commutative
period is simply the class $[\Omega]$ of the symplectic form $\Omega
\in \Omega^2(X/S)$, and it is easy to check that the non-commutative
period map is compatible with the base change.
\endproof

Comparing Theorem~\ref{main} and Theorem~\ref{kv}, we see that if
the quantization $\D_S$ of the universal deformation $\X/S$ is
canonical, then {\em all} quantizations of the symplectic manifold
$X/\Spec k$ can be obtained in a unique way by pullback from
$\D_S$. Thus $Q(X,\Omega_0)$ is identified with the set of sections
$$
s = \id \times s':\Delta \to \overline{S} = \Delta \times S
$$ 
of the canonical projection $\overline{S} \to \Delta$. The canonical
quantization of $X$ corresponds to the constant section $s = \id
\times [\Omega_0]$. Analogously, by Theorem~\ref{kv} every
symplectic deformation $X'/\Delta$ can be obtained by a pullback
with respect to a map $s:\Delta \to \overline{S}$ of the type
$$
s = \{0\} \times s':\Delta \to \overline{S} = \Delta \times S.
$$
In the case when the quantization $\D_S$ is not canonical,
Lemma~\ref{basech} no longer holds. However, it is easy to check
that the basic picture is still the same: any quantization of $X$
can be obtained in a unique way by pullback from $\D_S$. Thus even
in the case when $\X/S$ does not admit a canonical quantization, it
is still possible to fit all the quantizations into a single
multi-parameter deformation $\D_S$.

\section{Discussion.}\label{dis}

\renewcommand{\addcontentsline}[3]{\relax}

\subsection{} The main difference between our approach and that of
Fedosov is that Fedosov works in $C^\infty$ setting, where all
principal bundles with respect to nilpotent groups are
trivial. Therefore the group part of our Harish-Chandra torsor
$\M_q$ reduces to the bundle of symplectic frames in
$\T(X)$. Fedosov does this reduction implicitly, by choosing a
symplectic connection on $X$. After that, the only non-trivial part
of the quantization procedure is the construction of the flat
$\h$-valued connection. This can be done directly in the jet bundle
$J^\infty(D)$ associated to $D$. The principal bundle $\M_q$ itself
does not enter into the picture at any point.

\subsection{} Nest and Tsygan \cite{NT} generalize the Fedosov
construction to the holomorphic setting. In this case the principal
bundle can be quite non-trivial. However, it still reduces to the
symplectic frame bundle in the $C^\infty$ category. Nest and Tsygan
analyze the holomorphic non-triviality by considering the Dolbeault
complexes, and encoding everything into the $(0,1)$-part of the
Fedosov connection. Again, everything is done in the jet bundle, and
the principal bundle does not appear explicitly at any point.

\subsection{} De Wilde and Lecomte use a different approach -- they
choose an open cover of the manifold $X$ and glue together local
quantizations by explicit group-valued C\^{e}ch cocycles. Moreover,
Deligne re-tells their construction in the language of nonabelian
cohomology and gerbes. This is further from our approach in that it
does not use the jet bundles, but it is closer in that it does
mention groups and torsors more explicitly. The main difference is
in how to set up the induction process -- in other words, how to
filter the deformation problem by the powers of the Planck constant
$h$. Since De Wilde and Lecomte do not use jet bundles, they cannot
work directly with groups -- for all the groups and Lie algebras
that appear, they have to find an explicit description as
automorphisms and derivations of this or that algebraic
object. Unfortunately, it seems that it is not obvious how to
interpret our groups $(Der D)_p$ in this way. In particular, setting
$D_p = D/h^{p+1}$ gives an embedding
$$
(Der D)_p \to Der D_p,
$$ 
but this embedding is {\em not} an isomorphism (the difference
appears in the $h^p$ part -- the left-hand side contains only the
Hamiltonian vector fields there, while the right-hand side contain
all vector fields). This can be corrected at the first step by
imposing an additional bracket operation on $D_1$, the way we do in
Lemma~\ref{D1}. However, listing axioms for this bracket in higher
orders and constructing the deformation theory for such objects
seems to be quite cumbersome.  De Wilde-Lecomte and Deligne also use
some additional algebraic data, and this is partially successful, at
least in the $C^\infty$ setting where they work. But it does
introduce some complications into the proofs, and there are some
extra parasitic obstructions which one has to kill by hand.

\subsection{} The central role played by the central extension
\eqref{cent} is fully realized both by Fedosov and by De
Wilde-Lecomte. In De Wilde-Lecomte (retold by Deligne), it is used
to add necessary rigidity to $D_p$. In Fedosov, and even more so in
Nest-Tsygan, it appears in the connections themselves -- in our
notation, they are not $(\Der D)$-valued but $G$-valued. To
compensate for this, the connections are allowed to have non-trivial
curvature with values in the center of the Lie algebra $G$. It is
this curvature that parametrizes the quantizations. In our approach,
this appears coupled with the possible group-theoretic obstructions
in the guise of our non-commutative period map.

\subsection{} In general, the quantizations we construct are purely
formal. However, among manifolds admissible in our sense, one finds
compact smooth porjective manifolds over $\C$. In this situation, it
would be very interesting to try to use compactness and obtain some
sort of quantization which is analytic in $h$ in some appropriate
sense. We would like to note, though, that brute force does not
work: it is not possible to obtain a deformation of the sheaf
$\calo_X$ of holomorphic functions which is defined over an actual
small disc with coordinate $h$. Indeed, for every small open disc $U
\subset X$, the power series in $h$ which define the quantized
product of holomorphic functions on $U$ do converge. However, by
looking at the Weyl algebra it is elementary to check that the
radius of convergence roughly coincides with the size of
$U$. Therefore it goes to $0$ when $U$ is shrunk to a point.

\subsection{} One final word concerns a more explicit description of
the set $Q(X,\Omega)$ of isomorphism classes of quantizations. We
have embedded it canonically into $H^2_{DR}(X)[[h]]$, and we have
proved that $Q(X,\Omega)$ is non-canonically isomorphic to
$H^2_F(X)$. This is really weak -- essentially, we just say that
$Q(X,\Omega) \subset H^2_{DR}(X)[[h]]$ is a smooth algebraic
subvariety with correct trans\-versality properties
w.r.t. $H^2_F(X)[[h]] \subset H^2_{DR}(X)[[h]]$, and use the
implicit function theorem to get an identification $Q(X,\Omega)
\cong H^2_F(X)[[h]]$. In the commutative symplectic case considered
in \cite{KV}, the final answer is analogous (if the universal
symplectic deformation $\X/S$ admits a canonical quantization. the
answer is in fact literally the same, see
Section~\ref{ex}). However, at least for projective $X$ the full
answer is also known, due to the pioneering work of F. Bogomolov
\cite{bogo}. The period domain for commutative deformations of
irreducible holomorphic symplectic manifolds is a globally, not
infinitesemally defined quadric in $H^2_{DR}(X)$. This is a deep
result; in particular, we get a non-trivial and completely canonical
quadratic form on $H^2_{DR}(X)$, known as the Bogomolov-Beauville
form. What happens for non-commutative deformations? Nest and Tsygan
asked the same question, in their language. Moreover, they were able
to compute the ``first-order'' part of $Q(X,\Omega) \subset
H^2_{DR}(X)[[h]]$. The answer is expressed in terms of the so-called
{\em Rozansky-Witten} characteristic classes of the symplectic
manifold $X$; the reader can find it (without proof) in
Remark~\ref{rw}. It would be very interesting to obtain a full
answer. This would probably involve some non-linear combinations of
Rozansky-Witten classes -- hopefully no more than quadratic, but
possibly not.

\bigskip

\noindent
{\sc Northwestern University\\
Evanston, IL, USA\\
\mbox{}\\
Steklov Math Institute\\
Moscow, USSR}

\bigskip

\noindent
{\em E-mail addresses\/}: {\tt bezrukav@math.northwestern.edu}\\
\phantom{{\em E-mail addresses\/}: }{\tt kaledin@mccme.ru}

\end{document}